\documentclass{article}
\usepackage{graphicx}
\DeclareGraphicsExtensions{.jpg,.pdf}
\usepackage{amssymb,amsmath}
\usepackage[english]{babel}
\usepackage{color}

\begin{document}


\title{Well posedness of fluid-solid mixture models for biofilm spread}

\author{Ana Carpio (Universidad Complutense de Madrid),  \\
Gema Duro (Universidad Aut\'onoma de Madrid)}
\maketitle

{\bf Abstract}
Two phase solid-fluid mixture models are ubiquitous in biological
applications. For instance, models for growth of tissues and biofilms
combine time dependent and quasi-stationary boundary value problems 
set in domains whose boundary moves in response to variations in 
the mechano-chemical variables.
For a model of biofilm spread, we show how to obtain better posed 
models by characterizing the time derivatives of relevant quasi-stationary 
magnitudes in terms of additional boundary value problems. We also 
give conditions for well posedness of time dependent submodels set in 
moving domains depending on the motion of the boundary.  After 
constructing solutions for  transport, diffusion and elliptic submodels 
for volume fractions, displacements, velocities, pressures and 
concentrations with the required regularity, we  are able to handle the 
full model of biofilm spread in moving domains assuming we know the 
dynamics of the boundary. These techniques  are general and can be
applied in models with a similar  structure arising in biological 
and chemical engineering applications. \\

{\bf Keyword.}
Fluid-solid mixture models, thin film approximations, evolution
equations in moving domains, quasi-stationary approximations,
stationary transport equations. \\


\section{Introduction}  
\label{sec:intro}

Biofilms are bacterial aggregates that adhere to moist surfaces. Bacteria are
encased in a self-produced polymeric matrix \cite{biofilm} which shelters 
them from chemical and mechanical aggressions.
Biofilms formed on medical equipment, such as implants and catheters,
are responsible for hospital-acquired infections  \cite{hai}.
In industrial environments, they cause substantial economical and technical 
problems, associated to food poisoning, biofouling, biocorrosion,
contaminated ventilation  systems, and so on \cite{econ,slimy,zhu}. Modeling 
biofilm spread is important to be able to eradicate them.

\begin{figure}
\centering
\includegraphics[width=12cm]{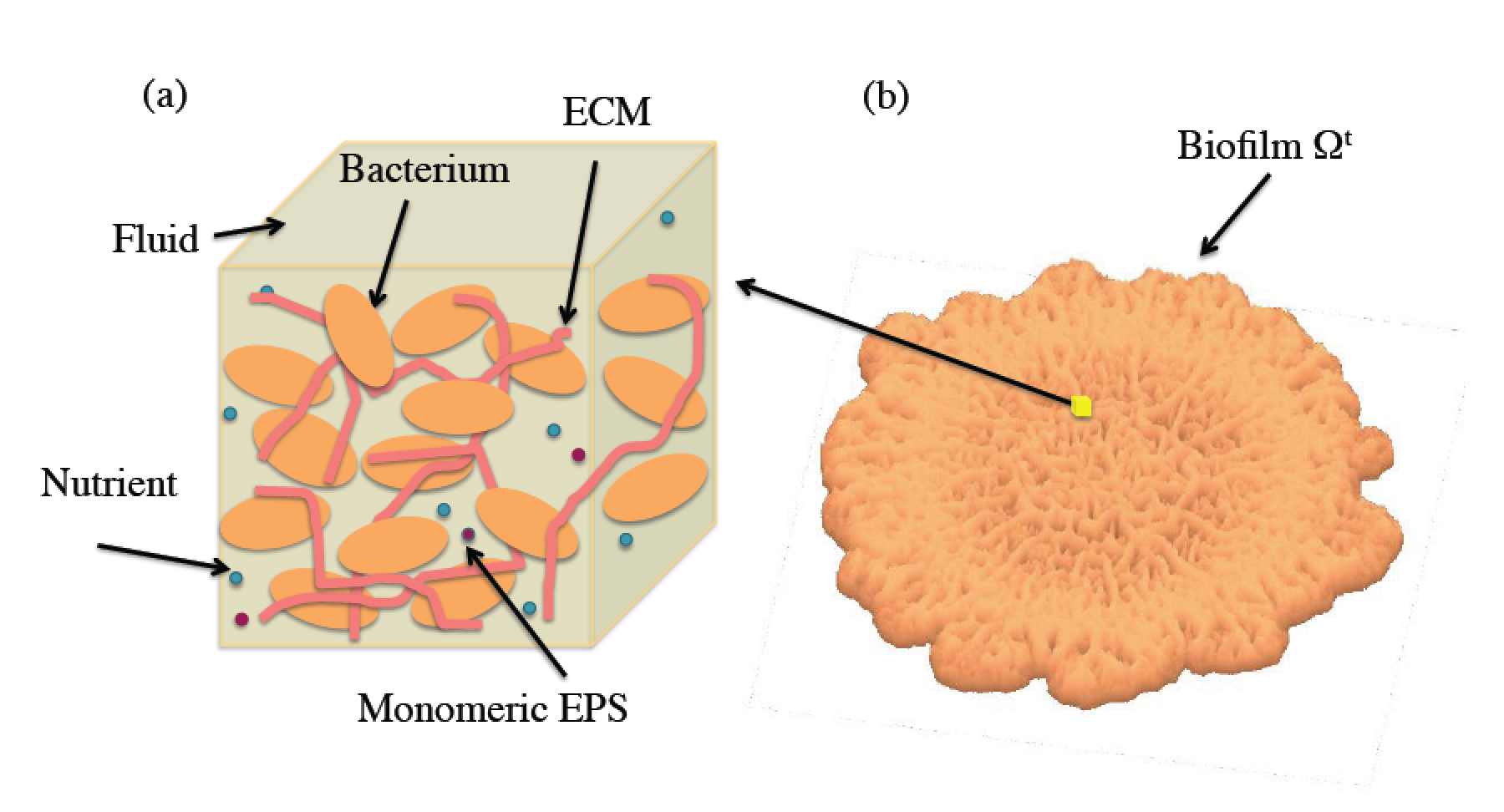}
\caption{(a) Schematic view of biofilm microscopic structure. Cells are
embedded in a network of polymeric threads forming the extracellular
matrix (ECM), while a liquid solution containing nutrients and chemicals
flows through the network.
(b) Schematic view of a biofilm spreading on a surface.}
\label{fig0}
\end{figure}

We describe here biofilms in terms of solid-fluid mixtures, see Figure
\ref{fig0}. At each point  $\mathbf x$ of the biofilm we have a solid fraction 
of biomass  $\phi_s(\mathbf x,t)$ (cell biomass, polymeric threads) and a 
volume  fraction of water $\phi_f(\mathbf x,t)$ containing dissolved substances 
(nutrients, autoinducers and so on), in such a way that $\phi_s(\mathbf x,t)
+\phi_f(\mathbf x,t)=1.$ The solid and fluid volume fractions move with 
velocities $\mathbf v_s$ and $\mathbf v_f$, respectively.

\begin{figure}
\centering
\includegraphics[trim=15mm 4.5cm 20mm 1cm,clip,width=6cm]{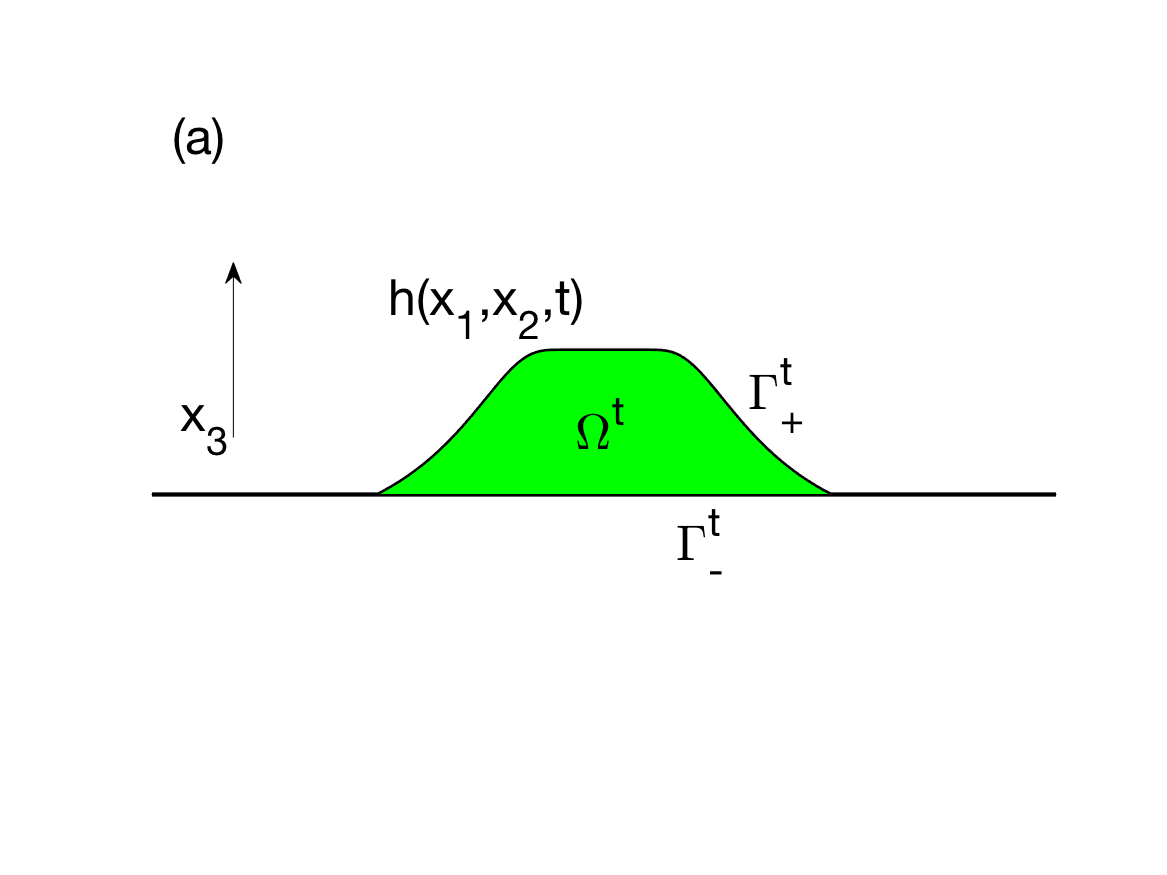}
\includegraphics[trim=15mm 4.5cm 20mm 1cm,clip,width=6cm]{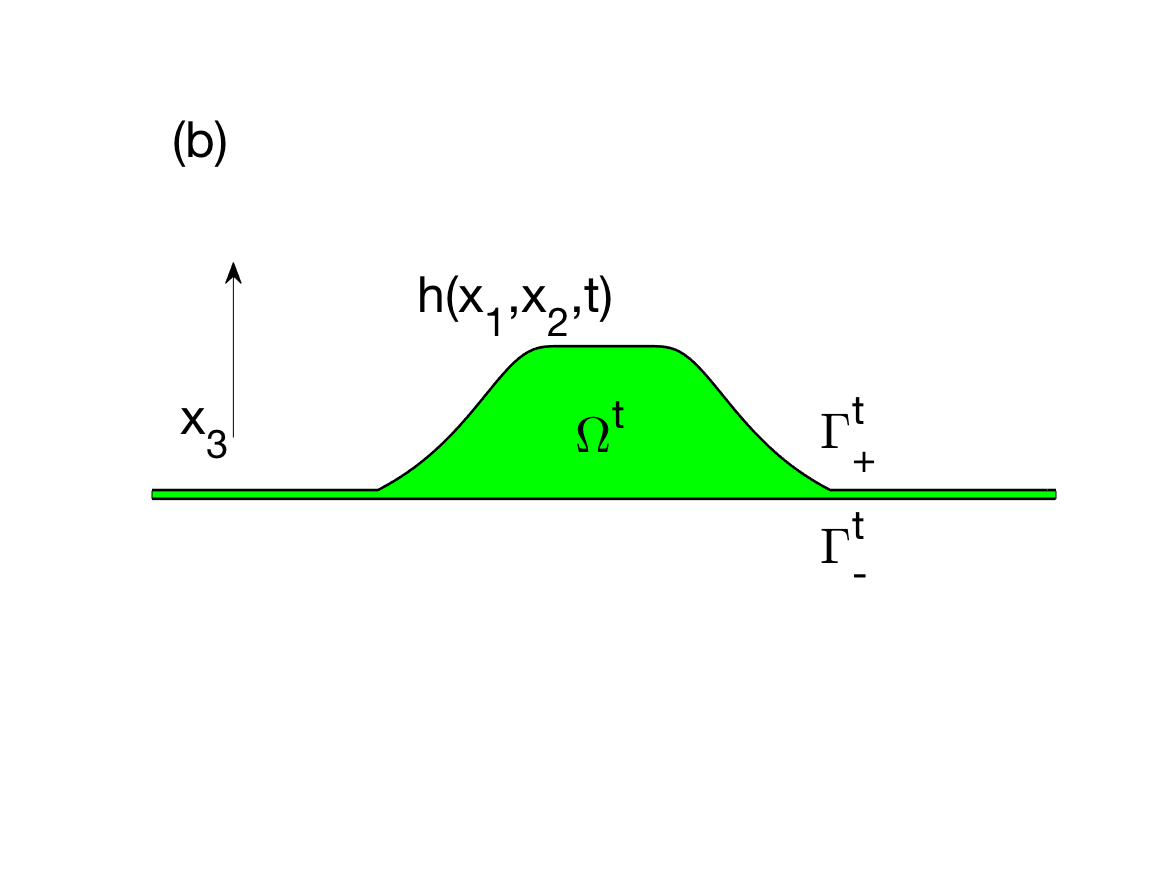}
\caption{Schematic representation of a biofilm slice $\Omega^t$ spreading 
on a surface (a) occupying a finite region and ending at triple contact points,
(b) spreading over precursor layers. The upper boundary $\Gamma^t_+$ 
represents the biofilm/air interface. The lower boundary $\Gamma^t_-$ 
represents the biofilm/agar interface, which provides nutrients and resources 
necessary for biofilm growth in our framework. }
\label{fig1}
\end{figure}

Biofilm spread on an air/solid interface is governed by the following  system 
of equations, see \cite{entropy,lanir}. Assume a biofilm occupies a region $\Omega^t$, 
that varies with time. Figure \ref{fig1} represents schematic views of two
dimensional slices. The upper boundary $\Gamma^t_+$ separates the biofilm 
from an outer fluid, that can be a liquid or air.  A lower boundary $\Gamma^t_-$ 
separates the biofilm  from the substratum it attaches to. The main variables satisfy
a set of quasi-stationary equations 
\begin{eqnarray} \begin{array}{lcl} \displaystyle
{\rm div}  ({\mathbf v}_f \phi_f) 
&=& - k_s {c \over c + K_s} \phi_s,  
\\ [1ex]  \displaystyle
{\rm div}(k_h(\phi_s) \nabla (p-\pi(\phi_s)) &=& 
{\rm div} (\mathbf v_s ),  
\\ [1ex]  \displaystyle
\mu \Delta \mathbf u_s + (\mu + \lambda) 
\nabla ({\rm div}(\mathbf u_s)) &=&  \nabla p, 
\\ [1ex]  \displaystyle
-d  \Delta c  + {\rm div}  (\mathbf v_f c)  
&=& - k_c { c \over c + K_c} \phi_s, 
\end{array} \label{ieqs}  
\end{eqnarray}
constrained by the additional conditions
\begin{eqnarray}
\phi_f \mathbf v_f  = - k_h(\phi_s) \nabla (p-\pi(\phi_s)) 
+ \phi_f \mathbf v_s, \quad
{\mathbf v_s} ={\partial \mathbf u_s \over \partial t},   
\quad \phi_f+\phi_s = 1, \label{iconds}
\end{eqnarray}
in the region occupied by the biofilm $\Omega^t$, which
varies with time.  
In this quasi-static framework, the displacement vector 
$\mathbf u_s(\mathbf x,t)$ and  the scalar pressure $p(\mathbf x,t)$,
volume fraction $\phi_s(\mathbf x,t)$ and concentration
$c(\mathbf x,t)$  fields depend on time through variations
of the  boundary $\Gamma^t$, which expands due to
cell division and swelling. The positive functions $k_h(\phi_s)$
and $\pi(\phi_s)$ represent the permeability and the osmotic 
pressure. This system is subject to a set of boundary conditions:
\begin{eqnarray} \begin{array}{ll} \displaystyle
p - \pi = p_{ext} - \pi_{ext},  & \qquad \mbox{on } \Gamma^t=
\Gamma^t_+\cup \Gamma^t_-,     
\\  \displaystyle
(\hat{\boldsymbol \sigma}(\mathbf u_s)  - p \mathbf I) \mathbf n 
= \mathbf t_{ext},  \; {\partial c \over \partial \mathbf n} =0,
& \qquad \mbox{on } \Gamma^t_+,   
\\ [1ex]  \displaystyle
\mathbf u_s = 0, \; c = c_0, & \qquad \mbox{on } \Gamma^t_-,    
\end{array}
 \label{ibcs}
\end{eqnarray}
where $\mathbf n$ is the outer unit normal and 
\[
\hat{\boldsymbol \sigma}(\mathbf u_s)=
\lambda {\rm Tr} 
(\boldsymbol  \varepsilon({\mathbf u}_s))\, \mathbf I +  
2 \mu \, \boldsymbol  \varepsilon({\mathbf u}_s),
\quad 
\varepsilon_{ij}({\mathbf u})= {1\over 2} \Big(  {\partial u_i \over \partial x_j } 
+ {\partial u_j \over \partial x_i}  \Big), \;  i,j=1,\ldots,n,
\]
$n=2, 3,$ represent elastic stress and strain tensors. Boundary 
conditions for $\phi_f$ are required or not depending on the sign 
of $\mathbf v_{f} \cdot \mathbf n$ at the border.
The displacement and velocity vectors have components
$\mathbf u = (u_1,\ldots,u_n)$ and $\mathbf v = (v_1,\ldots,v_n)$,
$n= 2,3$, respectively. All the parameters appearing in the model, 
$k_s$, $K_s$, $k_c$, $K_c$, $\mu$, $\lambda$,  $d$ are positive
constants. For ease of the reader, we have summarized
the modeling in Appendix A.  In some limits, the system can be
reformulated as a poroelastic model \cite{poroelastic, kapellos}.

The model is complemented with an equation for the dynamics
of $\Gamma^t$, $t>0$. If we consider biofilms represented
by the scheme in  Figure \ref{fig1}(a), the contact points between
biofilm, air and agar require specific additional  information to avoid
singularities. We will work with the geometry represented in
Figure \ref{fig1}(b), that avoids this difficulty by introducing precursor
layers \cite{seminara,degennes}. Then, $\Gamma^t_-$ is
fixed. The upper boundary  $\Gamma^t_+$  is parametrized
by a height function $h(x_1,x_2,t)$, which satisfies the equation
\cite{entropy}
\begin{eqnarray}
{\partial h \over \partial t } +
{\partial \over \partial x_1} \left[
\int_{0}^h (\mathbf v \cdot \hat{\mathbf x}_1) \, dx_3 \right]
+ {\partial \over \partial x_2} \left[
\int_{0}^h (\mathbf v \cdot \hat{\mathbf x}_2) \, dx_3 \right]  =
\mathbf v\cdot \hat{\mathbf x}_3\big|_0,
\label{iheight}
\end{eqnarray}
where the composite velocity of the mixture
$\mathbf v = \phi_f \mathbf v_f + \phi_s \mathbf v_s $ has 
components $\mathbf v \cdot \hat{\mathbf x}_i =
v_{s,i} - k_h(\phi_s)  {\partial (p-\pi) \over \partial x_i}, i=1,2,3$.

At present, only perturbation analyses and numerical studies  are
available for this type of models \cite{entropy, seminara} in simple
geometries. Asymptotic studies yield thin film type approximations 
for (\ref{ieqs})-(\ref{iheight}) assuming circular geometries and radial 
symmetry. Non standard lubrication  equations for the height $h$ 
are obtained, which admit families of self-similar solutions in radial 
geometries.  However, the construction of reliable  numerical solutions 
of the model in general experimental configurations faces difficulties
due to the lack of well-posedness results.

In this paper, we assume we know the dynamics of the upper boundary 
$\Gamma_+^t$, given by a smooth curve $x_3 = h(x_1,x_2,t)$, 
and develop an existence and stability theory for the model equations.
To simplify the analysis, we take $k_h(\phi_s) = k_h >0$,
$k_h(\phi_s)/\phi_f = \xi_\infty >0$ and $\pi(\phi_s) = \Pi \phi_s >0$. 
In this quasi-stationary framework, the displacements
$\mathbf u_s$ depend on time through the motion of the boundary.
However, we lack equations for the velocities, other than the
relation ${\partial \mathbf u_s \over \partial t} = \mathbf v_s$. In Section
\ref{sec:velocity} we obtain a system of equations characterizing the 
velocity:
\begin{eqnarray} 
{\rm div}(\hat{\boldsymbol \sigma}(\mathbf v_s)) =
\mu \Delta \mathbf v_s + (\mu + \lambda) 
\nabla ({\rm div}(\mathbf v_s)) =  \nabla p_t, 
& \; \mbox{in $\Omega^t$},  \nonumber \\  
\mathbf  v_s = 0,  & \; \mbox{on $\Gamma_-^t$},   
 \label{eivs} \\  
 \hat{\boldsymbol \sigma}(\mathbf v_s) \mathbf n 
= {\partial \mathbf g\over \partial t} +  \mathbf r(\mathbf g,\mathbf u_s), 
&\; \mbox{on $\Gamma_+^t$,}
\nonumber
\end{eqnarray}
with $g= - p \mathbf n  = -( p_{\rm ext} - \pi_{\rm ext})\mathbf n$
and  ${\mathbf r}$ to be defined later. A similar equation is obtained for 
$p_t$ from the equation for $p$.
Taking the divergence of the equations for $\mathbf u_s$ and
$\mathbf v_s$ we find additional equations to close the system
\begin{eqnarray}
{de \over dt} = k_h  (2 \mu + \lambda)  \Delta e -  k_h \Pi \Delta \phi_s, 
& \; \mbox{in $\Omega^t$}, 
\label{eie} \\
{de_t \over dt} = k_h (2 \mu + \lambda)  \Delta e_t -  k_h \Pi \Delta \phi_{s,t},
& \; \mbox{in $\Omega^t$}, 
\label{eiet}
\end{eqnarray}
where $e = {\rm div} (\mathbf u_s)$  and  $e_t = {\rm div} (\mathbf v_s)$. 
We will neglect $\Delta \phi_{s,t}$ in (\ref{eiet}) because $\Pi$ and 
$\Delta \phi_s$ are small compared to other terms. 
Notice that (\ref{eie}) and (\ref{eiet}) are time dependent problems set
in time dependent domains, while most results in the literature refer
to fixed domains.

The construction of solutions for such systems combines a number of 
difficulties that we will address in stages.
Section \ref{sec:velocity} characterizes the time derivatives of $\mathbf u_s$ 
and $p$,  solutions of elliptic problems in time dependent domains, by 
means of additional boundary value problems. In this way we improve the
stability of the model, since solving additional partial differential equations in 
each spatial domain is more effective than approximating time derivatives 
by quotients of differences of solutions calculated in variable spatial domains.
Section \ref{sec:tdsubmodels} establishes well posedness results
for  linear parabolic problems (\ref{eiet}) set in domains with moving 
boundaries for specific types of parametrizations. 
Section \ref{sec:stsubmodels} considers the elliptic and stationary transport 
problems involved in the quasi-stationary submodels, separately and in fixed 
domains, under hypotheses motivated by asymptotic studies and numerical
solutions. Finally, section \ref{sec:full} considers the full coupled time 
dependent problem  and section  \ref{sec:conclusions} discusses our
conclusions and open issues.
A final appendix summarizes modeling details.


\section{Differentiation of quasi-stationary problems}
\label{sec:velocity}

In the previous section, we have defined the velocity $\mathbf v_s$ as
the time derivative of the displacement $\mathbf u_s$. The change
in time of $\mathbf u_s$ is due to the motion of the upper boundary
$\Gamma^t_+$, that is, time variations in $h$. In this section
we seek an equation characterizing $\mathbf v_s$.
We expect $\mathbf v_s$ to solve the same boundary value problem 
as $\mathbf u_s$, but differentiating all sources with respect to time. 
However, since the boundary $\Gamma^t$ of $\Omega^t$ moves with 
time,  we need to calculate the adequate  boundary conditions too.

In the region $\Omega^t$ occupied by the moving biofilm, the displacements 
$\mathbf u_s$ of the solid phase satisfy equations (\ref{ieqs}) with boundary
conditions (\ref{ibcs}).
To simplify later computations, it is convenient to recast these 
equations in the general linear elasticity framework. The components 
of the displacement $u_{j}(t), j=1,\ldots,n$, $n$ being the dimension, 
fulfill
\begin{eqnarray} 
- {\partial \over \partial x_\alpha } \left(
c_{j \alpha m \beta} {\partial u_{m}(t) \over \partial x_{\beta}   } \right)
= f_j(t),  \quad j=1,\ldots,n, & \qquad\mbox{in $\Omega^t$}, 
\nonumber \\  
 u_j(t) = 0,  \quad j=1,\ldots,n, & \qquad\mbox{on $\Gamma_d^t$}, 
\label{elastostatic} \\  
\displaystyle  c_{j \alpha m \beta}   {\partial u_{m}(t)  \over
\partial x_{\beta} }  n_{\alpha}(t)  
=g_j(t), \quad j=1,\ldots,n, &\qquad\mbox{on $\Gamma_n^t$,}
\nonumber
\end{eqnarray}
where $\mathbf n(t)$ is the outer unit normal vector
and $c_{j \alpha m \beta}$ the elastic constants. 
$\Gamma_n^t$ and $\Gamma_d^t$ are parts of the boundary $\Gamma^t$ 
where we enforce conditions on the stresses of the displacements,
respectively.

We use the Einstein summation convention that implies summation 
over a set of indexed terms in a formula when repeated in it. In the
above equations, summation over $\alpha, \beta, m$ is implied, but not 
over $j$.
The elastic constants  $c_{j\alpha m \beta}$ for a isotropic solids
like the ones we consider are
\begin{equation*}
c_{j \alpha m \beta}=\lambda \delta_{j \alpha}\delta_{m \beta}+\mu
(\delta_{jm}\delta_{\alpha \beta} +\delta_{j \beta}\delta_{\alpha m})
\label{elasticisotropic}
\end{equation*}
where $\delta_{jm}$ stands for the Kronecker delta, whereas $\lambda$
and $\mu$ represent the Lam\'e constants. The stress tensor is
\begin{equation*}
\sigma_{j\alpha} = c_{j \alpha m \beta} \varepsilon_{m\beta}
= \lambda \delta_{j \alpha} \varepsilon_{pp} + 2 \mu \varepsilon_{j\alpha}.
\label{elasticstress}
\end{equation*}


In this framework, the velocity $\mathbf v$  is the `Fr\`echet derivative' 
or `domain derivative' of $\mathbf u$ with respect to $t$ \cite{dirichlet2D}, 
which is characterized by the solution of a boundary value problem, as we
show next.
\\

{\bf Theorem 2.1.} {\it We assume that the  body $\mathbf f$ and boundary 
$\mathbf g$ forces are differentiable in time, with values in $[L^2(\Omega^t)]^n$ 
and $[L^2(\Gamma^t)]^n$, respectively, with $t>0$, $n=2,3$  being the
dimension. Moreover, the $C^2$ boundaries $\Gamma^t$ are obtained deforming
$\Gamma^0$ along a smooth vector field $\boldsymbol \nu$. Then, the time derivative 
$\mathbf v(t)= {\partial \mathbf u(t) \over \partial t}$, $t>0$, of the displacement 
given by (\ref{elastostatic}) satisfies
\begin{eqnarray} 
- {\partial \over \partial x_\alpha} \left(
c_{j \alpha m \beta} {\partial v_{m}(t) \over \partial x_{\beta} } \right)
= {\partial f_j(t) \over \partial t},  
\quad j=1,..,n, & \quad \mathbf x \in \Omega^t, 
\nonumber \\  
 v_j(t) = 0,  \quad j=1,..,n, & \quad \mathbf x \in \Gamma_d^t,   
 \label{delastostatic} \\  
\displaystyle  c_{j \alpha m \beta}   {\partial v_{m}(t)  \over
\partial x_{\beta}}  n_{\alpha}(t) 
= {\partial g_j(t) \over \partial t} +  r_j(g_j(t),\mathbf u(t)), 
\quad j=1,..,n, &\quad \mathbf x \in \Gamma_n^t,
\nonumber
\end{eqnarray}
where 
\begin{eqnarray} \begin{array}{l}
\displaystyle
r_j =c_{j\alpha m\beta} {\partial u_m(t)\over \partial x_\beta} 
{\partial \nu_q \over \partial x_\alpha} n_q(t)
+ c_{j\alpha m\beta}  {\partial u_m(t) \over \partial x_\beta} 
{\partial (\nu_p   n_\alpha(t)) \over \partial x_p}  
\\[1ex] \displaystyle
 \hskip 5mm + c_{j \alpha m \beta} {\partial u_m(t) \over \partial x_{\beta} }
{\partial \nu_p \over \partial x_p} n_\alpha(t) 	 
- g_j(t) \mathbf n(t)^T \nabla \boldsymbol \nu \, \mathbf n(t),
\quad j=1,\ldots,n. \end{array} \label{rj}
\end{eqnarray}
}


As a corollary, we get the expressions of interest for our model. \\

{\bf Corollary 2.2.} {\it Under the previous hypotheses, 
the time derivative $\mathbf v_s(t)$, $t>0$, of the solution 
$\mathbf u_s$ of (\ref{ieqs}) with boundary conditions
(\ref{ibcs}) satisfies
\begin{eqnarray} 
{\rm div}( \hat{\boldsymbol \sigma}(\mathbf v_s))=
\mu \Delta \mathbf v_s + (\mu + \lambda) 
\nabla ({\rm div}(\mathbf v_s)) =  \nabla p_t, 
& \; \mathbf x \in \Omega^t,  \nonumber \\  
\mathbf  v_s = 0,  & \; \mathbf x \in \Gamma_-^t,   
 \label{evs} \\  
 \hat{\boldsymbol \sigma}(\mathbf v_s) \mathbf n 
= {\partial g_j \over \partial t} +  r_j(g_j,\mathbf u_s), 
\; j=1,2 &\; \mathbf x \in \Gamma_+^t, 
\nonumber
\end{eqnarray}
with $g= - p \mathbf n  = -( p_{\rm ext} - \pi_{\rm ext})\mathbf n$
and  ${\mathbf r}$ is defined by {(\ref{rj})} with 
$c_{j \alpha m \beta}=
\lambda \delta_{j \alpha}\delta_{m \beta}+\mu(\delta_{jm}\delta_{\alpha \beta} 
+\delta_{j \beta}\delta_{\alpha m}).$}\\

{\bf Corollary 2.3} {\it Under the previous hypotheses, 
assuming $k_k(\phi_s)=k_h$ and $\pi(\phi_s) = \Pi \phi_s$,
the derivative $p_t(t)= {\partial p(t) \over \partial t}$, 
$t>0,$ of the solution  $p$ of (\ref{ieqs}) with Dirichlet boundary
conditions $p= p_{ext}(t)$ satisfies,
\begin{eqnarray} \begin{array}{ll}
 k_h \Delta  p_t  =  {\rm div}(\mathbf v_{s,t})
 +  k_h \Pi \Delta  \phi_{s,t} +  ,  & \mathbf x \in \Omega^t, 
\nonumber \\  [1ex]
p_t  = p_{\rm ext}'(t),   & \mathbf x \in \Gamma^t. 
\end{array} \label{pt} 
\end{eqnarray}} 

{\bf Proof of Theorem 2.1.}
We will follow a similar variational approach to that employed in 
\cite{dirichlet2D} for 2D exterior elasticity problems with zero Dirichlet 
boundary conditions on a moving boundary.  We are going to calculate
the derivative at $t=0$. Similar arguments hold for any $t>0$.

{\it Step 1: Variational formulation.} First, we write the boundary value
problem for $\mathbf u$ in variational form \cite{raviart}.
The boundary value problem (\ref{elastostatic}) becomes: Find
$\mathbf u^t \in [H^1_{\Gamma_d^t}(\Omega^t)]^n$ such that
\begin{equation}\label{Variationale}
\displaystyle  b^t(\Omega^t; \mathbf u^t, \mathbf w^t)= {\ell}^t(\Omega^t;\mathbf w^t),
\qquad \forall \, \mathbf w^t \in [H^1_{\Gamma_d^t}(\Omega^t)]^n,
\end{equation}
where
\begin{eqnarray}
b^t(\Omega^t; \mathbf u^t,\mathbf w^t) = \int_{\Omega^t}  c_{j \alpha m \beta}
{\partial u_{m}^t\over \partial x_{\beta}^t } {\partial
\overline w_{j}^t \over \partial x_{\alpha}^t  }   \, d\mathbf x^t, 
\qquad  \forall \, \mathbf u^t,\,\mathbf w^t \in [H^1_{\Gamma_d^t}(\Omega^t)]^n, 
\label{bt} \\
\ell^t(\Omega^t; \mathbf w^t)= \int_{\Omega^t}  f_j(t) w_j^t \, d\mathbf x^t
+ \int_{\Gamma_n^t}  g_j(t) w_j^t \, d\mathbf S_{\mathbf x^t}, 
 \qquad\forall \, \mathbf w^t \in [H^1_{\Gamma_d^t}(\Omega^t)]^n.
\label{ellt}
\end{eqnarray}
Here, $H^1_{\Gamma_d^t}(\Omega^t)$ denotes the usual Sobolev space of 
$H^1(\Omega^t)$ functions vanishing on $\Gamma_d^t \subset \partial \Omega^t$. 
$H^1(\Omega^t)$ if formed by all functions whose square, and the squares
of their derivatives, are integrable in $\Omega^t$, that is, belong to $L^2(\Omega^t).$
When $\mathbf f(t) \in [L^2(\Omega^t)]^n$, $\mathbf g \in [L^2(\Gamma^t)]^n$ and 
${\rm meas}(\Gamma_d^t)\neq 0$, this problem admits a unique solution 
$\mathbf u^t \in [H^1_{\Gamma_d^t}(\Omega^t)]^n$ \cite{raviart}, 
which in fact belongs to $[H^2(\Omega^t)]^n$, vanishes on ${\Gamma_d^t}$ 
and satisfies
$\boldsymbol \sigma (\mathbf u^t) \, \mathbf n = \mathbf g $ on  $\Gamma_n^t= 
\partial \Omega^t \setminus \Gamma_d^t.$ For $t=0$, we have $u^0.$
Here, $\sigma_{\alpha j} (\mathbf u^t) = c_{j \alpha m \beta} {\partial u_{m}^t \over \partial x_{\beta}}.$

{\it Step 2: Change of variables.} 
We now transform all the quantities appearing in (\ref{bt})-(\ref{ellt}) back to 
the initial configuration $\Omega^0$. The process is similar
to transforming deformed configurations  back to a reference configuration in 
continuum mechanics \cite{gurtin}. 
We are assumig that the evolution of the moving part of the
boundary $\Gamma^t  = \left\{ \mathbf x + t \, \boldsymbol 
\nu(\mathbf x) \, | \, \mathbf x \in \Gamma^0  \right\}$
is given by a family of deformations  $\mathbf x^t = \phi^t(\mathbf x)
= \mathbf x + t \, \boldsymbol \nu(\mathbf x)$
starting from a smooth surface $ \Gamma^0 \in C^2$ (twice
differentiable) and following a smooth vector field
$\boldsymbol \nu \in C^2 (\Omega)$, $\Omega^t \subset \Omega$,
$t>0$. 
The  deformation gradient is the jacobian of the change of variables 
\cite{feijoooberai}
\begin{eqnarray}
\mathbf J^t(\mathbf x) = \nabla_{\mathbf x} \boldsymbol \phi^t(\mathbf x)  = 
\left({\partial x^t_i \over \partial x_j}(\mathbf x) \right) =
\mathbf I + t \, \nabla \boldsymbol \nu(\mathbf x),
\label{jacobian}
\end{eqnarray} 
and its inverse $ (\mathbf J^t)^{-1}   
= \left({\partial x_i \over \partial x^t_j}\right)$ is the jacobian of the inverse
change of variables.
Then, volume and surface elements are related by
\begin{eqnarray}
d \mathbf x^t = {\rm det} \, \mathbf J^t(\mathbf x) \, d \mathbf x, \quad
d S_{\mathbf x^t} = {\rm det} \, \mathbf J^t(\mathbf x)
\|\, (\mathbf J^t(\mathbf x))^{-T} \mathbf n \| dS_{\mathbf x}
\label{changevolume}
\end{eqnarray}
and the chain rule for derivatives reads 
$\nabla_{\mathbf x} u_m(\mathbf x^t(\mathbf x)) = 
\mathbf (J^t(\mathbf x))^T \nabla_{\mathbf x^t} u_m(\mathbf x^t(\mathbf x)),$
that is,
$\nabla_{\mathbf x^t} u_m = (\mathbf J^t)^{-T} \nabla_{\mathbf x} u_m$.
For each component we have
\begin{eqnarray}
{\partial u_m \over \partial x_\beta^t}(\mathbf x^t(\mathbf x))  = 
{\partial u_m \over \partial x_k}(\mathbf x^t(\mathbf x)) 
(J^t)^{-1}_{k\beta}(\mathbf x).
\label{changederivative}
\end{eqnarray}

%

We define $\tilde {\mathbf u}(\mathbf x)= \mathbf u^t  \circ \phi^t (\mathbf x)
= \mathbf u^t (\mathbf x^t(\mathbf x))$, definition that extends to 
$\tilde {\mathbf w}$ and  other functions. Changing variables and using (\ref{changevolume})-(\ref{changederivative}) we have:
\begin{eqnarray}
b^t(\Omega^t; \mathbf u^t,\mathbf w^t)   
 = \int_{\Omega^t}  c_{j \alpha m \beta}
{\partial u_m^t\over \partial x_{\beta}^t }(\mathbf x^t) {\partial
  w_j^t \over \partial x_{\alpha}^t  }(\mathbf x^t)   \, d\mathbf x^t = 
 \nonumber  \\
\int_{\Omega^0}  c_{j \alpha m \beta}
{\partial \tilde u_m\over \partial x_{p}  }(\mathbf x) 
(J^t)^{-1}_{p \beta}(\mathbf x)
{\partial \tilde w_j  \over \partial x_{q} }(\mathbf x) 
(J^t)^{-1}_{q \alpha}(\mathbf x)  \, 
{\rm det} \, \mathbf J^t(\mathbf x) \, d \mathbf x =
\tilde b^t(\Omega^0; \tilde{\mathbf u},\tilde{\mathbf w})
 \label{changeb} \\
\ell^t(\Omega^t; \mathbf w^t)  = 
\int_{\Omega^t}  f_j(\mathbf x^t,t) w_j^t(\mathbf x^t) \, d\mathbf x^t  +
\int_{\Gamma_n^t}  g_j(\mathbf x^t,t) w_j^t(\mathbf x^t) 
\, d S_{\mathbf x^t} = \nonumber  \\ 
 \int_{\Omega^0}  \hskip -2mm \tilde f_j(\mathbf x,t) \tilde w_j(\mathbf x) \, 
{\rm det} \, \mathbf J^t \, d \mathbf x  \!+\!
\int_{\Gamma_n^0}  \hskip -2mm \tilde g_j(\mathbf x,t) \tilde w_j(\mathbf x)
{\rm det} \, \mathbf J^t
\|\, (\mathbf J^t)^{-T} \! \mathbf n \| dS_{\mathbf x}
= \tilde \ell^t(\Omega^0; \tilde{\mathbf w}).
\label{changeell}
\end{eqnarray}
For arbitrary test functions $\mathbf w^t \in [H^1_{\Gamma_d^t}(\Omega^t)]^n$,
$\tilde {\mathbf w} \in [H^1_{\Gamma_d^t}(\Omega^0)]^n$ is a test function in
$\Omega^0$. Therefore, we obtain the equivalent variational formulation:
Find $\tilde{\mathbf u} \in [H^1_{\Gamma_d^t}(\Omega^0)]^n$ such that
\begin{equation}\label{Variationalehat}
\displaystyle  \tilde b^t(\Omega^0; \tilde{\mathbf u}, \mathbf w)=
\tilde {\ell}^t(\Omega^0;\mathbf w),
\qquad \forall \, \mathbf w \in [H^1_{\Gamma_d^t}(\Omega^0)]^n,
\end{equation}
with $ \tilde b^t(\Omega^0; \tilde{\mathbf u}, \mathbf w)$ and 
$\tilde {\ell}^t(\Omega^0;\mathbf w)$ defined in (\ref{changeb})-(\ref{changeell})
replacing $\tilde{\mathbf w}$ by $\mathbf w$.

Let us analyze the dependence on $t$ of the terms appearing in the
expression for $\tilde b^t$ and $\tilde \ell^t$. From the definitions of the Jacobian 
matrices (\ref{jacobian}) we obtain \cite{feijoooberai,dirichlet2D}
\begin{eqnarray}
{\rm det} \, \mathbf J^t(\mathbf x)  = 1 + t \, {\rm div}(\boldsymbol \nu(\mathbf x) ) + O(t^2),  
\label{expanddet}\\
(\mathbf J^t)^{-1}(\mathbf x)  = \mathbf I - t \, \nabla \boldsymbol \nu(\mathbf x)  + O(t^2), 
\label{expandinv} \\
{\rm det} \, \mathbf J^t(\mathbf x) \, \| (\mathbf J^t(\mathbf x))^{-T} \mathbf n \| =
1 + t \, {\rm div}_{\Gamma}(\boldsymbol \nu(\mathbf x)) + O(t^2),
\label{expandsurf}
\end{eqnarray}
where ${\rm div}_{\Gamma}(\boldsymbol \nu(\mathbf x)) =
{\rm div}(\boldsymbol \nu(\mathbf x)) - \mathbf n^T \nabla \boldsymbol \nu(\mathbf x)
\mathbf n.$
Inserting (\ref{expanddet})-(\ref{expandinv}) in (\ref{changeb}) we find the following
expansions. When $p=\beta$ and $q=\alpha$ we get
\begin{eqnarray}
\int_{\Omega^0}  c_{j \alpha m \beta}
{\partial \tilde u_m\over \partial x_{\beta}  }
{\partial   w_j  \over \partial x_{\alpha} } \, d \mathbf x
+ t  \int_{\Omega^0}  c_{j \alpha m \beta}
{\partial \tilde u_m \over \partial x_{\beta}  }
{\partial w_j \over \partial x_{\alpha} }
\, {\rm div}(\boldsymbol \nu)  \, d \mathbf x 
\label{expandequal}\\
- t \int_{\Omega^0}  c_{j \alpha m \beta} \left[
{\partial \tilde u_m\over \partial x_{\beta}  } 
{\partial \nu_\beta \over \partial x_\beta} 
{\partial w_j  \over \partial x_{\alpha} } 
+
{\partial \tilde u_m \over \partial x_{\beta}  } 
{\partial w_j  \over \partial x_{\alpha} } 
{\partial \nu_\alpha \over \partial x_\alpha} \right] \, d \mathbf x
+ O(t^2), \nonumber 
\end{eqnarray}
whose leading term is $b^0(\Omega^0; \tilde{\mathbf u}, \mathbf w ) $.
When $p\neq \beta$ and $q \neq \alpha$ the summands
are $O(t^2)$. The remaining terms provide the  contribution
\begin{eqnarray*}
-t \int_{\Omega^0}  c_{j \alpha m \beta} \left[ 
{\partial \tilde u_m \over \partial x_{p}  }
{\partial \nu_{p} \over \partial x_\beta}
{\partial w_j  \over \partial x_{\alpha} }
+
{\partial \tilde u_m \over \partial x_\beta  }
{\partial w_j  \over \partial x_{q} }
{\partial \nu_{q} \over \partial x_\alpha} \right]
 \, d \mathbf x
 + O(t^2),
\end{eqnarray*}
with $p \neq \beta$, $q=\alpha$ in the first one and
$q \neq \alpha$, $p=\beta$ in the second one.
Adding up the contributions we get 
\begin{eqnarray}
\tilde b^t(\Omega^0; \tilde{\mathbf u}, \mathbf w )   =
b^0(\Omega^0; \tilde{\mathbf u}, \mathbf w ) +
t[I_1(\tilde{\mathbf u})+I_2(\tilde{\mathbf u})+I_3(\tilde{\mathbf u})]
+O(t^2),
\label{expandb}
\end{eqnarray}
where
\begin{eqnarray} 
\begin{array}{ll}
I_1(\tilde{\mathbf u}) =& \int_{\Omega^0}  c_{j \alpha m \beta}
{\partial \tilde u_m \over \partial x_\beta  } 
{\partial w_j  \over \partial x_\alpha } 
\, {\rm div}(\boldsymbol \nu)  \, d \mathbf x,  \\
I_2(\tilde{\mathbf u})  =& - \int_{\Omega^0}  c_{j \alpha m \beta}  
{\partial \tilde u_m \over \partial x_{p}  } 
{\partial \nu_{p} \over \partial x_\beta} 
{\partial w_j  \over \partial x_\alpha } \, d \mathbf x,  \\
I_3(\tilde{\mathbf u}) =&  
-  \int_{\Omega^0}  
c_{j \alpha m \beta}  {\partial \tilde u_m \over \partial x_\beta} 
 {\partial w_j  \over \partial x_q } 
{\partial \nu_q \over \partial x_\alpha} \, d \mathbf x  =
\int_{\Omega^0}  {\partial \over \partial x_\alpha}
\left(c_{j \alpha m \beta}  {\partial \tilde u_m \over \partial x_\beta} 
\right) {\partial w_j  \over \partial x_q } 
 \nu_q \, d \mathbf x \\ 
 & + \int_{\Omega^0}   c_{j \alpha m \beta}  
{\partial \tilde u_m\over \partial x_\beta} 
{\partial^2 w_j  \over \partial x_\alpha \partial x_q} 
 \nu_q \, d \mathbf x 
  -  \int_{\partial \Omega^0}  c_{j \alpha m \beta}  
{\partial \tilde u_m \over \partial x_{\beta}  } 
n_\alpha  {\partial w_j \over \partial x_q } 
 \nu_q  \, d S_\mathbf x.
 \end{array} \label{expandIs}
\end{eqnarray}
Similarly, from the definition (\ref{changeell}) of the linear form 
$\tilde \ell^t$ and the definition of 'material derivative' $\dot{\mathbf f}$
\begin{eqnarray}
\tilde{\mathbf f}(\mathbf x,t) = {\mathbf f}(\mathbf x^t(\mathbf x),t)
= \mathbf f(\mathbf x,0)  + t \,
\dot{\mathbf f}(\mathbf x,0) + O(t^2),
\label{materialf}
\end{eqnarray}
we find the expansion
\begin{eqnarray}
\begin{array}{lll}
\tilde \ell^t(\Omega^0;  \mathbf w )  
&= \int_{\Omega^0}  f_j(0) w_j \, d \mathbf x
+ t \int_{\Omega^0}  [ f_j(0) {\rm div}(\boldsymbol \nu)  +
\dot f_j(0) ] w_j \, d \mathbf x  
 \\ [1.5ex]
&+ \int_{\Gamma_n^0}   g_j(0) w_j \, d S_{\mathbf x} + 
t \int_{\Gamma_n^0}  [ g_j(0) {\rm div}_\Gamma(\boldsymbol \nu)  
+ \dot g_j(0) ] w_j \, d S_{\mathbf x} + O(t^2)
\end{array}
\label{expandell}
\end{eqnarray}
whose leading term is $\ell^0(\Omega^0; \mathbf w).$

{\it Step 3. Variational problem for the domain derivative $\mathbf u'$.}
Let us compare the transformed function $\tilde{\mathbf u} $ and the solution 
$\mathbf u^0$ of 
$b^0(\Omega^0;\mathbf u^0, \mathbf w) = \ell^0(\Omega^0;\mathbf w)$.
For any $\mathbf w \in [H^1_{\Gamma_d^t}(\Omega^0)]^n$ we have
\begin{eqnarray}
b^0(\Omega^0;\tilde{\mathbf u}- \mathbf u^0, \mathbf w) =
b^0(\Omega^0;\tilde{\mathbf u}, \mathbf w) - \ell^0(\Omega^0;\mathbf w) = 
\nonumber \\
b^0(\Omega^0;\tilde{\mathbf u}, \mathbf w) 
- \tilde b^t(\Omega^0;\tilde {\mathbf u}, \mathbf w)
+ \tilde \ell^t(\Omega^0;\mathbf w) - \ell^0(\Omega^0;\mathbf w).
\label{expandvariational}
\end{eqnarray}
Well posedness of the variational problems (\ref{Variationale}) with respect to 
changes in domains $\Omega^t$ and sources $\mathbf f(t), \mathbf g(t)$, 
implies uniform bounds on the solutions for $t \in [0,T]$:
$\|\mathbf u^t\|_{[H^1(\Omega^t)]^n} \leq C(T)$, 
$\| \tilde{\mathbf u} \| _{[H^1(\Omega^0)]^n} \leq C(T).$
Expansions (\ref{expandb})-(\ref{expandell}) show that the  right hand side
in (\ref{expandvariational})  tends to zero as $t \rightarrow 0$.
Well posedness of the variational problem again implies $\tilde{\mathbf u}
\rightarrow \mathbf u^0$ in $[H^1_{\Gamma_d^t}(\Omega^0)]^n$
as $t\rightarrow 0$.

Dividing by $t$ equation (\ref{expandvariational}) and using
(\ref{expandb})-(\ref{expandell}), we find
\begin{eqnarray}
\begin{array}{lll}
b^0(\Omega^0;{\tilde{\mathbf u}- \mathbf u^0 \over t}, \mathbf w) &= &
{1\over t} [b^0(\Omega^0;\tilde{\mathbf u}, \mathbf w) - 
\tilde b^t(\Omega^0;\tilde {\mathbf u}, \mathbf w)] + 
{1\over t} [ \tilde \ell^t(\Omega^0;\mathbf w) -\ell^0(\Omega^0;\mathbf w)]   
\\ [1.5ex]
&= &  - [I_1(\tilde{\mathbf u})+I_2(\tilde{\mathbf u})+I_3(\tilde{\mathbf u})]+
\int_{\Omega^0}  [ f_j(0) {\rm div}(\boldsymbol \nu)  + \dot f_j(0) ]   w_j 
\, d \mathbf x  
\\ [1.5ex]
&+ &  \int_{\Gamma_n^0}  [ g_j(0) {\rm div}_\Gamma(\boldsymbol \nu)  +
\dot g_j(0) ] w_j \, d S_{\mathbf x}+ O(t).
\end{array} \nonumber
\end{eqnarray}
Then, the limit 
$\dot {\mathbf u} = {\rm lim}_{t \rightarrow 0} {\tilde{\mathbf u}- 
\mathbf u^0 \over t}$ satisfies
\begin{eqnarray} 
\begin{array}{ll}
b^0(\Omega^0;\dot {\mathbf u}, \mathbf w) = &
\int_{\Omega^0} 
[ f_j(0) {\rm div}(\boldsymbol \nu)  + \dot f_j(0) ]   w_j   d \mathbf x
- [I_1(\mathbf u^0)\!+\!I_2(\mathbf u^0)\!+\!I_3(\mathbf u^0)] 
\\ [1.5ex]
& + \int_{\Gamma_n^0}  [ g_j(0) {\rm div}_\Gamma(\boldsymbol \nu)  
+ \dot g_j(0) ] w_j \, d S_{\mathbf x}.
\end{array}
\label{equdot} 
\end{eqnarray}
As before, the function $\dot {\mathbf u}$ is the so called `material derivative',
that is, $\dot {\mathbf u} = {\partial \mathbf u \over \partial t} + 
\nabla \mathbf u^0 \, \boldsymbol \nu$.  The domain derivative becomes
$\mathbf u' = \dot {\mathbf u} - \nabla \mathbf u^0 \, \boldsymbol \nu$. 
 Then,\begin{eqnarray}
b^0(\Omega^0; \mathbf u', \mathbf w) =
b^0(\Omega^0;\dot {\mathbf u}, \mathbf w)
- b^0(\Omega^0;\nabla \mathbf u^0 \, \boldsymbol \nu, \mathbf w),
\label{equprime}
\end{eqnarray}
where 
\begin{eqnarray*}
b^0(\Omega^0;\nabla \mathbf u^0 \, \boldsymbol \nu, \mathbf w) =
\int_{\Omega^0}  {\partial \over \partial x_\beta}
\left( c_{j \alpha m \beta}  {\partial u_m^0\over \partial x_p }  \nu_{p} \right)
{\partial w_j \over \partial x_\alpha}  \, d \mathbf x.
\end{eqnarray*}
Notice that this function vanishes on $\Gamma_d$ whenever $\dot {\mathbf u}$ 
and $\boldsymbol \nu$ do so.

{\it Step 4. Differential equation for the domain derivative ${\mathbf u}'$.}
We evaluate the different terms in the right hand side of (\ref{equdot})
to calculate the right hand side in (\ref{equprime}). First, notice that
$-{\partial \over \partial x_\alpha} \left(c_{j \alpha m \beta}  
{\partial  u_m^0 \over \partial x_{\beta}} \right) = f_j(0) $ in $\Omega^0$ and 
$ c_{j \alpha m \beta}  {\partial u_m^0 \over \partial x_{\beta}  } 
n_\alpha= g_j(0)$ on $\Gamma_n^0$, $u_j^0=0$ on $\Gamma_d^0$,
$j=1,...,n$, imply:
\begin{eqnarray*}
\begin{array}{ll}
I_3(\mathbf u^0) 
 & = \int_{\Omega^0}  \left( {\partial f_j(0) \over \partial x_q}  \nu_q + f_j(0)
{\partial \nu_q \over \partial x_q} \right)
w_j    \, d \mathbf x 
- \int_{\partial \Omega^0}  f_j(0)   w_j  n_q \nu_q   \, d \mathbf x \\ [1.5ex]
 & -  \int_{\Gamma_n^0}  g_j(0)  {\partial w_j \over \partial x_q } 
 \nu_q  \, d S_\mathbf x
+ \int_{\Omega^0}   c_{j \alpha m \beta}  
{\partial u_m^0\over \partial x_\beta} 
{\partial^2 w_j  \over \partial x_q  \partial x_\alpha} 
 \nu_q  \, d \mathbf x. 
\end{array} \nonumber
\end{eqnarray*}  
Using
${\partial  u_m^0\over \partial x_p} {\partial \nu_p \over \partial x_\beta} =  
{\partial \over \partial x_\beta} \left({\partial  u_m^0 \over \partial x_p } 
 \nu_p \right) - {\partial^2  u_m^0\over \partial x_p \partial x_\beta} \nu_p $,
we get
\begin{eqnarray*}
\begin{array}{ll}
I_2(\mathbf u^0) 
& = - b^0(\Omega^0;\nabla \mathbf u^0 \, \boldsymbol \nu, \mathbf w) 
 - \int_{\Omega^0}   c_{j \alpha m \beta} 
{\partial u_m^0\over \partial x_{\beta} }   {\partial \nu_{p} 
\over \partial x_p} {\partial  w_j \over \partial x_\alpha}  \, d \mathbf x  
\\ [1.5ex]
& - \int_{\Omega^0}   c_{j \alpha m \beta} 
{\partial  u_m^0\over  \partial x_\beta }    \nu_{p} 
{\partial^2 w_j  \over \partial x_\alpha \partial x_p} \, d \mathbf x  
+   \int_{\partial \Omega0}   c_{j \alpha m \beta} 
{\partial  u_m^0\over  \partial x_\beta }  \nu_{p} 
n_p {\partial w_j  \over \partial x_\alpha} \, d \mathbf x.
\end{array} \nonumber
\end{eqnarray*}
As a result of the two previous identities
\begin{eqnarray*}
\begin{array}{ll}
- [I_1(\mathbf u^0)+I_2(\mathbf u^0)+I_3(\mathbf u^0)]= 
b^0(\Omega^0;\nabla \mathbf u^0 \, \boldsymbol \nu, \mathbf w) 
-    \int_{\partial \Omega^0}     c_{j \alpha m \beta} 
{\partial  u_m^0\over  \partial x_\beta }  \nu_{p} 
n_p {\partial w_j  \over \partial x_\alpha} \, d \mathbf x   \\[1.5ex]
-\int_{\Omega^0}    \left( {\partial f_j(0) \over \partial x_q}  \nu_q 
+ f_j(0) {\partial \nu_q \over \partial x_q} \right) w_j    \, d \mathbf x 
+ \int_{\partial \Omega^0}   f_j(0)   w_j  n_q \nu_q   \, d \mathbf x
+  \int_{\Gamma_n^0}  g_j(0)  {\partial w_j \over \partial x_q } 
 \nu_q  \, d S_\mathbf x 
\end{array} \nonumber 
\end{eqnarray*}
and (\ref{equprime}) becomes
\begin{eqnarray*}
\begin{array}{l}
b^0(\Omega^0; \mathbf u', \mathbf w) =
-    \int_{\partial \Omega^0}   c_{j \alpha m \beta} 
{\partial  u_m^0\over  \partial x_\beta }  \nu_{p} 
n_p {\partial w_j  \over \partial x_\alpha} \, d \mathbf x 
+ \int_{\partial \Omega^0}  f_j(0)   w_j  n_q \nu_q   \, d \mathbf x +
\\  [1.5ex]
\int_{\Omega^0}  f_j'(0)   w_j  \, d \mathbf x + 
\int_{\Gamma_n^0}  [ g_j(0) {\rm div}_\Gamma(\boldsymbol \nu)  
+ \dot g_j(0) ] w_j \, d S_{\mathbf x}
+  \int_{\Gamma_n^0}  g_j(0)  {\partial w_j \over \partial x_q } 
 \nu_q  \, d S_\mathbf x.
\end{array}
\end{eqnarray*}
Integrating by parts in  $b^0(\Omega^0; \mathbf u', \mathbf w)$
and choosing $\mathbf w$ with compact support inside $ \Omega^0$, 
this identity yields  the following equation for $\mathbf u'$  in $\Omega^0$
\begin{eqnarray}
- {\partial \over \partial x_\alpha}
\left(c_{j \alpha m \beta} {\partial u_{m}'  \over   
\partial x_{\beta}}(\mathbf x) \right) =  
f_j'(\mathbf x,0), \quad j=1,...,n.
\label{equprimediff}
\end{eqnarray}
However, to obtain a pointwise boundary condition for $\mathbf u'$
we need to rewrite the integral on $\partial \Omega^0$ in such a way 
that no derivatives of the test function $\mathbf w$ are involved.

{\it Step 5: Boundary condition for the domain derivative $\mathbf u'$.}
We integrate by parts the original expressions of $I_i(\mathbf u^0)$,
$i=1,2,3$ to get
\begin{eqnarray*} 
\begin{array}{lll}
I_1 & = &
- \int_{\Omega^0}  {\partial \over \partial x_\alpha}\left( c_{j \alpha m \beta}
{\partial   u_m^0\over \partial x_{\beta}  } \, {\rm div}(\boldsymbol \nu)  \right)
w_j \, d \mathbf x  
+ \int_{\partial \Omega^0}   c_{j \alpha m \beta}
{\partial u_m^0 \over \partial x_{\beta} }
\, {\rm div}(\boldsymbol \nu) n_\alpha w_j \, d S_\mathbf x,
\end{array}  
\end{eqnarray*} \vskip -5mm
\begin{eqnarray*}
\begin{array}{lll}
I_2  &= &   
- \int_{\Omega^0} c_{j\alpha m\beta} {\partial \over \partial x_\beta}
\left({\partial u_m^0 \over \partial x_p}  \nu_p  \right)
{\partial w_j \over \partial x_\alpha} d \mathbf x  
- \int_{\Omega^0} c_{j\alpha m\beta} {\partial \over \partial x_\alpha}
{\partial \over \partial x_p}
\left( {\partial u_m^0 \over \partial x_\beta} \nu_p  \right)  
w_j   d \mathbf x \\ [1.5ex]
& + & \int_{\Omega^0} c_{j\alpha m\beta} {\partial \over \partial x_\alpha}
\left( {\partial u_m^0 \over \partial x_\beta} {\partial \nu_p  \over
\partial x_p} \right)  w_j   d \mathbf x
+ \int_{\partial \Omega^0} c_{j\alpha m\beta} 
{\partial^2 u_m^0 \over \partial x_p \partial x_\beta}  \nu_p  
w_j n_\alpha d \mathbf x 
\end{array} 
\end{eqnarray*} \vskip -4mm
\begin{eqnarray*}
\begin{array}{lll}
I_3 &= &
\int_{\Omega^0} {\partial \over \partial x_q} {\partial \over \partial
x_\alpha} \left( c_{j\alpha m\beta} {\partial u_m^0 \over \partial x_\beta} 
  \nu_q  \right)  w_j  d \mathbf x 
+ \int_{\Omega^0} {\partial \over \partial x_q} \left(
f_j(0) \nu_q  \right)  w_j   d \mathbf x \\ [1.5ex]
&  - & \int_{\partial \Omega^0} 
c_{j\alpha m\beta} {\partial u_m^0 \over \partial x_\beta} 
{\partial \nu_q \over \partial x_\alpha} n_q w_j  d S_{\mathbf x} 
\end{array}
\end{eqnarray*}

Adding up to compute $ -[I_1+I_2+I_3]$,  integrating by parts 
$b^0(\Omega^0;  \mathbf u', \mathbf w)$, inserting (\ref{equprimediff})  
in (\ref{equprime}) and setting $\boldsymbol \nu=0$ on 
$\Gamma_d$ we find
\begin{eqnarray*}
\begin{array}{ll}
\int_{\Gamma_n^0} c_{j \alpha m \beta}{\partial u_m' \over \partial x_\beta}      
 n_\alpha w_j \, d S_{\mathbf x} &=  \int_{\partial \Omega^0} 
c_{j\alpha m\beta} {\partial u_m^0 \over \partial x_\beta} 
{\partial \nu_q \over \partial x_\alpha} n_q w_j  \, d S_{\mathbf x} \\ [1.5ex]
& - \int_{\partial \Omega^0} c_{j\alpha m\beta}  [
{\partial^2 u_m^0 \over \partial x_p \partial x_\beta}  \nu_p  
n_\alpha +  {\partial u_m^0 \over \partial x_{\beta} }
\, {\rm div}(\boldsymbol \nu) n_\alpha]  w_j \, d S_\mathbf x \\ [1.5ex]
& + \int_{\Gamma_n^0}  [ g_j(0) {\rm div}_\Gamma(\boldsymbol \nu)  
+ \dot g_j(0) ] w_j \, d S_{\mathbf x}.
\end{array}
\end{eqnarray*}
Now, using identifies $
- c_{j\alpha m\beta}  {\partial^2 u_m^0 \over \partial x_p \partial x_\beta}  
\nu_p  n_\alpha = - {\partial \over \partial x_p} (g_j(0) \nu_p) 
+ c_{j\alpha m\beta}  {\partial u_m^0 \over \partial x_\beta} 
{\partial (\nu_p   n_\alpha) \over \partial x_p}, $ and
$  g_j(0) {\rm div}_\Gamma(\boldsymbol \nu)  
+ \dot g_j(0)  = {\partial \over \partial x_p} (g_j(0) \nu_p) - g_j(0)
\mathbf n^T \nabla \boldsymbol \nu \, \mathbf n + g'_j(0), $
we obtain
\begin{eqnarray} 
\begin{array}{lll} 
c_{j \alpha m \beta}  {\partial   u_m' \over \partial x_{\beta}} n_\alpha &  = &
c_{j\alpha m\beta} {\partial u_m^0 \over \partial x_\beta} 
{\partial \nu_q \over \partial x_\alpha} n_q
+ c_{j\alpha m\beta}  {\partial u_m^0 \over \partial x_\beta} 
{\partial (\nu_p   n_\alpha) \over \partial x_p}  \\ [1.5ex]
& & + c_{j \alpha m \beta} {\partial u_m^0 \over \partial x_{\beta} }
{\partial \nu_p \over \partial x_p} n_\alpha 	 
- g_j(0) \mathbf n^T \nabla \boldsymbol \nu \, \mathbf n + g'_j(0)
\end{array} \label{bcprime}
 \end{eqnarray}
on $\Gamma_n^0.$ $\square$


\section{Study of  diffusion problems in time dependent domains}
\label{sec:tdsubmodels}

We study here parabolic problems of the form 
\begin{eqnarray} \begin{array}{lll}
e_t - \kappa \Delta e = f(\mathbf x, t), & \; \mathbf x \in \Omega^t, & t>0, \\
e = g(t), & \; \mathbf x \in \Gamma^t, & t>0, \\
e(\mathbf x, 0) = e_0, &  \; \mathbf x \in \Omega^t. &
\end{array} \label{tdmodel}
\end{eqnarray}


As in Section \ref{sec:velocity} we assume that the evolution of the moving 
part of the boundary is given by a family of deformations  \cite{feijoooberai}
\begin{eqnarray*}
\Gamma^t  = \left\{ \mathbf x + t \, \boldsymbol \nu(\mathbf x) \, | \, 
\mathbf x \in \Gamma^0  \right\},
\end{eqnarray*}
starting from a smooth surface $ \Gamma^0  \in C^2$ (twice differentiable) and 
following a smooth vector field $\boldsymbol \nu \in C(\overline{\Omega^0}) \cup C^2(\Omega^0)$. 
We can assume $e(t)=0$ by making the change $e = \hat e + g$.
Then $\hat e$ solves (\ref{tdmodel}) with zero Dirichlet boundary condition,
initial datum $e_0(\mathbf x)-g(0)$ and right hand side $f(\mathbf x, t)-g'(t)$.
Therefore, we will work with zero Dirichlet boundary conditions in the sequel.
To solve (\ref{tdmodel}) we will first refer it to a fixed domain and then
construct  converging Faedo-Galerkin approximations.

\subsection{Variational formulation in the undeformed configuration}
\label{sec:undeformed}

As usual, we denote as $H^1_{0}(\Omega^t)$ the subspace of 
$H^1(\Omega^t)$ formed by functions whose trace vanishes on 
$\Gamma^t$ with the induced norm.
Multiplying (\ref{tdmodel}) by $w^t \in H^1_{0}(\Omega^t)$ and 
integrating, we find
\begin{eqnarray*}
\begin{array}{l} \begin{array}{l} \displaystyle
\int_{\Omega^t} e_t(\mathbf x^t,t) w^t(\mathbf x^t) \, d \mathbf x^t \!+\!
\int_{\Omega^t} \nabla_{\mathbf x^t} e(\mathbf x^t,t) \nabla_{\mathbf x^t}
w^t(\mathbf x^t)  \, d \mathbf x^t \!=\! 
\int_{\Omega^t} f(\mathbf x^t,t) w^t(\mathbf x^t) \, d \mathbf x^t
\end{array} \end{array} \label{var_t}
\end{eqnarray*}
for each $t$. We use (\ref{jacobian}), (\ref{changevolume}), 
(\ref{changederivative}) to refer these integrals to a fixed domain.

We define $\tilde {w}(\mathbf x)= w^t  \circ \phi^t (\mathbf x)
= w^t (\mathbf x^t(\mathbf x))$, $\phi^t $ as in (\ref{jacobian}).
Notice that
\begin{eqnarray} \begin{array}{l} \displaystyle
e_t(\mathbf x^t(\mathbf x),t) = {d\over dt}[e(\mathbf x^t(\mathbf x),t)] 
- \nabla_{\mathbf x^t} e(\mathbf x^t(\mathbf x),t)^T {d \mathbf x^t \over dt} 
\\[2ex] \displaystyle
= {d\over dt} \tilde e(\mathbf x,t) - (\mathbf J^t)^{-T}
\nabla_{\mathbf x} \tilde e(\mathbf x,t)^T \tilde{\boldsymbol \nu}(\mathbf x).
\end{array}  \label{changederivativet} 
\end{eqnarray}
After changing variables, problem (\ref{tdmodel}) reads: Find
$e \in C([0,T],L^2(\Omega^0)) \cap L^2(0,T;H^1_0(\Omega^0))$ 
such that $e(\mathbf x, 0) = e_0(\mathbf x)$ and
\begin{eqnarray} \begin{array}{l} \displaystyle
\int_{\Omega^0} \tilde e_t(\mathbf x,t) \tilde w(\mathbf x) \, 
{\rm det} \, \mathbf J^t(\mathbf x) \, d \mathbf x 
- \int_{\Omega^0} \nabla_{\mathbf x} \tilde e(\mathbf x,t)^T
(\mathbf J^t(\mathbf x))^{-1} \tilde{\boldsymbol \nu}(\mathbf x)
\tilde w(\mathbf x) \, {\rm det} \, \mathbf J^t(\mathbf x) \, d \mathbf x 
\\ [2ex] \displaystyle
+ \int_{\Omega^0} \nabla_{\mathbf x} e(\mathbf x,t)^T
((\mathbf J^t(\mathbf x))^T \mathbf J^t(\mathbf x))^{-1}
\nabla_{\mathbf x} e(\mathbf x,t) \tilde w(\mathbf x) \,
 \,  {\rm det} \, \mathbf J^t(\mathbf x) \, d \mathbf x
\\ [2ex] \displaystyle
= \int_{\Omega^0} \tilde f(\mathbf x,t) \tilde w(\mathbf x) \, 
{\rm det} \, \mathbf J^t(\mathbf x) \, d \mathbf x.
\end{array} \label{var_0}
\end{eqnarray}
Since $w^t \in H^1_{0}(\Omega^t)$, we have 
$\tilde w \in H^1_{0}(\Omega^0)$. In fact, we can take the same
arbitrary function $w \in H^1_{0}(\Omega^0)$ for all $t$.


\subsection{Construction of stable solutions}
\label{sec:tdexistence}

Consider a basis  $\{\phi_1, \phi_2, \ldots, \phi_M \ldots \}$ of the  
Hilbert space $L^2(\Omega)$. We  choose the normalized
eigenfunctions $\phi_j \in H^2(\Omega) \cap H^1_{0}(\Omega)$,
$j\in \mathbb N$, of $-\Delta$ in $H^1_{0}(\Omega)$, see  \cite{brezis}.
\vskip 1mm

{\bf Theorem 3.1} {\it Let $\Omega \subset \mathbb R^n$ be an open and 
bounded $C^2$ domain. Given a function $f \in C^1([0,T]; L^2(\Omega))$ 
there exists a unique solution
$u \in C([0,T]; H^2(\Omega)) \cap H^1(0,T;H^1_{0}(\Omega))$ of
\begin{eqnarray} \begin{array}{l} \displaystyle
\int_{\Omega}  u_t(\mathbf x,t)  w(\mathbf x) \, c(\mathbf x,t) \, d \mathbf x 
+ \int_{\Omega}  \nabla u(\mathbf x,t)^T \mathbf b(\mathbf x,t) w(\mathbf x)
 \, d \mathbf x + \\ [2ex] \displaystyle
\int_{\Omega} \nabla u(\mathbf x,t)^T 
\mathbf A(\mathbf x,t) \nabla w(\mathbf x) \, d \mathbf x   
= \int_{\Omega}  f(\mathbf x,t)   w(\mathbf x)   \, d \mathbf x,
\end{array} \label{var_ref}
\end{eqnarray}
for all $w \in H^1_{0}(\Omega)$, $t \in [0,T]$, provided 
\begin{itemize}
\item 
$\mathbf A(\mathbf  x,t) \in C^1(\overline{\Omega} \times [0,T])$, 
$\mathbf b(\mathbf x,t) \in C^1(\overline{\Omega} \times [0,T])$ and
$c(\mathbf x,t) \in C^2(\overline{\Omega} \times [0,T])$,
\item the matrices $\mathbf C^M(t)$ with elements 
$\int_\Omega  c(t) \phi_m \phi_k  d \mathbf x$, $m,k=1, \ldots, M,$ 
are invertible for $t \in [0,T]$, 
\item  the matrices $\mathbf A(\mathbf x, t)$ are uniformly coercive, that is,
$\boldsymbol \xi^T  \mathbf A(\mathbf x, t) \boldsymbol \xi \geq a_0 
|\boldsymbol \xi|^2,$
$a_0>0,$ for all $\boldsymbol \xi \in \mathbb R^n$, and the scalar field
$c(\mathbf  x,t)$ is bounded from below, $c(\mathbf  x,t) \geq c_0 >0$,
for all $\mathbf x \in \mathbb R^n$ and $t>0$,
\item $u_0 \in L^2(\Omega)$ and 
$w_0 ={\rm div}(\mathbf A(\mathbf x, 0) \nabla u_0(\mathbf x))
+ \mathbf b(\mathbf x, 0)^T \nabla u_0(\mathbf x) \in L^2(\Omega)$. 
\end{itemize}
Moreover, the solution depends continuously on parameters
and data.  } \\
We obtain a solution for the original time dependent problem set in 
a moving domain undoing the change of variables.
\vskip 1mm

{\bf Proof.} {\it Existence.}
We use the Faedo-Galerkin method \cite{lionsmagenes,lions}.
First, we change  variables $u(\mathbf x,t) = e^{\lambda t} v(\mathbf x,t)$,
$u_t(\mathbf x,t) = e^{\lambda t} [v_t(\mathbf x,t) + \lambda v(\mathbf x,t)]$,
with $\lambda >0$ to be selected large enough. We obtain  similar variational
equations for $v$ with an additional term $\lambda c v$ and $g$ and $f$ multiplied
by $e^{-\lambda t}$. Then we seek approximate solutions
$v^M(\mathbf x,t) = \sum_{m=1}^M \alpha_m(t) \phi_m(\mathbf x)$ 
such that
\begin{eqnarray} 
\begin{array}{l} \displaystyle
 \int_{\Omega} c(\mathbf x, t) \,v^M_t(\mathbf x,t)  
 w(\mathbf x)  \, d \mathbf x 
  + \int_{\Omega} 
\sum_{p,q=1}^n a_{pq}(\mathbf x,t)
{\partial v^M \over \partial x_p}(\mathbf x,t) 
{\partial w \over \partial x_q}(\mathbf x)  \, d \mathbf x 
 \\[2ex] \displaystyle
 + \int_{\Omega} 
\lambda  c(\mathbf x,t) v^M(\mathbf x,t) w(\mathbf x,t) \, d \mathbf x 
 + \int_{\Omega}  \sum_{p=1}^n b_p(\mathbf x,t)
{\partial v^M \over \partial x_p}(\mathbf x,t)    w(\mathbf x) 
\, d \mathbf x 
 \\[2ex] \displaystyle
 = \int_{\Omega}  e^{-\lambda t} f(\mathbf x,t)  w(\mathbf x) \,  
d \mathbf x, \\[2ex]
\displaystyle 
v^M(\mathbf x, 0) =  \sum_{m=1}^M \alpha_m(0) \phi_m(\mathbf x), \quad 
\alpha_m(0) = \int_\Omega u_0(\mathbf x) \phi_m(\mathbf x) d \mathbf x,
\end{array}  \label{var_x}
\end{eqnarray}
for all $w \in V^M={\rm span}\{\phi_1, \phi_2, \ldots, \phi_M\}$.
We find a system of $M$ differential equations for the coefficient functions 
$\alpha_m(t)$ setting $w = \phi_k$, $k=1,\ldots,M$,
\begin{eqnarray} \begin{array}{r} \displaystyle
\sum_{m=1}^M \alpha_m'(t) \int_\Omega c(t) \phi_m \phi_k \, d \mathbf x 
=  - \sum_{m=1}^M \alpha_m(t)  \int_{\Omega}  \sum_{p=1}^n b_p(t)
 {\phi_m\over \partial x_p} \phi_k \, d \mathbf x  - 
\\[2ex] \displaystyle
\sum_{m=1}^M \alpha_m(t) \int_\Omega \left[ \sum_{p,q=1}^n a_{pq}(t)
{\partial \phi_m \over \partial x_p} {\partial \phi_k \over \partial x_q} 
+ \lambda c(t) \phi_m  \phi_k \right]  \, d \mathbf x 
+ \int_{\Omega}  e^{-\lambda t} f(t) \phi_k \,  d \mathbf x.
\end{array}  \label{var_xm}
\end{eqnarray}
This can be written as a linear system with continuous and bounded 
coefficients in $[0,T]$
\begin{eqnarray*} 
 {d \over dt} \boldsymbol \alpha^{M} 
= \mathbf C^M(t)^{-1} \mathbf A^M(t) \boldsymbol \alpha^M  +  
\mathbf C^M(t)^{-1} \mathbf g^M(t)  + \mathbf C^M(t)^{-1}\mathbf f^M(t) 
\end{eqnarray*}
with initial datum $\boldsymbol \alpha^M(0)$, which admits a unique
solution $\boldsymbol \alpha^M(t)$, $t \in [0,T]$ \cite{coddington}.
Multiplying  identity (\ref{var_xm}) by $\alpha_k$ and adding over $k$, 
we obtain
\begin{eqnarray} 
\begin{array}{l} \displaystyle
{1\over 2}{d\over dt} \int_{\Omega} c(\mathbf x, t) |v^M(\mathbf x,t)|^2  
d \mathbf x 
 + \int_{\Omega} \sum_{p,q=1}^n a_{pq}(\mathbf x,t)
{\partial v^M \over \partial x_p}(\mathbf x,t) 
{\partial v^M \over \partial x_q}(\mathbf x,t)  \, d \mathbf x 
\\[2ex]  \displaystyle 
\!+\! \int_{\Omega}  \hskip -1mm  (\lambda c - {1\over 2}  c_t)(\mathbf x, t) 
|v^M(\mathbf x,t)|^2  \, d \mathbf x 
+ \int_{\Omega}  \sum_{p=1}^n b_p(\mathbf x,t)
{\partial v^M \over \partial x_p}(\mathbf x,t) v^M(\mathbf x,t)
\, d \mathbf x 
  \\[2ex]
\displaystyle = 
\int_{\Omega}  e^{-\lambda t} f(\mathbf x,t)  v^M(\mathbf x,t) \,  
d \mathbf x.  
\end{array}  \label{var_vM}
\end{eqnarray}
Integrating in $[0,T]$ and using coercivity, lower bounds for $A$ and
$c$, $L^\infty$ bounds, as well as Young's inequality \cite{brezis}, we find
\begin{eqnarray*}
{c_0\over 2} \int_{\Omega} |v^M(\mathbf x,t)|^2  d \mathbf x
+ {a_0 \over 2} \int_0^t \int_{\Omega} |\nabla v^M(\mathbf x,t)|^2  d \mathbf x ds
+ {\lambda c_0 \over 2} \int_0^t \int_{\Omega} |v^M(\mathbf x,t)|^2  d \mathbf x ds
\\ 
\leq  {\|c \|_{L^\infty_{xt}}\over 2}\ \|v^M(0)\|_{L^2(\Omega)} 
+ {1\over 2}\ \| f\|_{L^2(0,T,L^2(\Omega))}
\end{eqnarray*}
for $\lambda$ large enough depending on $a_0$, $\|c_t \|_{L^\infty_{xt}}$, 
$c_0$, $\|\mathbf b\|_{L^\infty_{xt}}$, $n$.
Gronwall inequality, and the fact that $v^M(0) \rightarrow u_0$ in 
$L^2$, imply  that $v^M$ is bounded in $L^\infty(0,T,L^2(\Omega))$ and 
$L^2(0,T,H^1_{0}(\Omega))$. We extract a subsequence $v^{M'}$ 
converging a limit $v$ weakly star in $L^\infty(0,T,L^2(\Omega))$ and 
weakly in $L^2(0,T,H^1_{0}(\Omega))$.
Moreover, ${d\over dt}\int_\Omega c(t) v^{M'}(t) \phi_k \, d \mathbf x$ tends 
to ${d\over dt}\int_\Omega c(t) v(t) \phi_k \, d \mathbf x$ in the sense of 
distributions in ${\cal D}'(0,T)$ for any $k$. Similar convegences hold for
$u^{M'}$ and $u= e^{\lambda t}v$. We undo the change in (\ref{var_x}), 
multiply by a function $\psi \in C_c^\infty([0,T))$, integrate over $t$ and 
pass to the limit as $M' \rightarrow \infty$ to find
\begin{eqnarray*} 
\begin{array}{r} \displaystyle
-  \int_{\Omega} c(\mathbf x, 0) \, u(\mathbf x,0)  
 w(\mathbf x) \psi(0) \, d \mathbf x 
- \int_0^t \int_\Omega c_t(\mathbf x,t) u(\mathbf x,t) w(\mathbf x) 
\psi(t) \, d \mathbf x ds +
 \\[2ex] \displaystyle
 \int_0^t \int_{\Omega} \left[ \sum_{p,q=1}^n a_{pq}(\mathbf x,t)
{\partial v \over \partial x_p}(\mathbf x,t) 
{\partial w \over \partial x_q}(\mathbf x) +
\sum_{p=1}^n b_p(\mathbf x,t)
{\partial v^M \over \partial x_p}(\mathbf x,t)    w(\mathbf x)\right]
\psi(t) \, d \mathbf x  ds
 \\[2ex] \displaystyle
 =  \int_0^t  \int_{\Omega}  e^{-\lambda t} f(\mathbf x,t)  w(\mathbf x) 
\psi(t) \,  d \mathbf x ds,
\end{array}  
\end{eqnarray*}
for any $w \in H^1_{0}(\Omega)$, so that the limiting solution 
satisfies the condition on the initial data and the equation
\begin{eqnarray}
c u_t - {\rm div}(\mathbf A \nabla u) + \mathbf b^T \nabla u =  f
\label{distributions}
\end{eqnarray}
in the sense of distributions \cite{lions,raviart}.

{\it Uniqueness.}
To prove uniqueness, we assume there are two solutions $u_1$ and
$u_2$, and set $u=u_1-u_2$. We subtract the equations satisfied by
both, multiply by $u$, set $u=e^{\lambda t}v$ and
integrate over $\Omega$ to get
\begin{eqnarray*} 
\begin{array}{l} \displaystyle
{1\over 2}{d\over dt} \int_{\Omega} c(\mathbf x, t) |v(\mathbf x,t)|^2  
d \mathbf x 
+ \int_{\Omega} \sum_{p,q=1}^n a_{pq}(\mathbf x,t)
{\partial v \over \partial x_p}(\mathbf x,t) 
{\partial v \over \partial x_q}(\mathbf x,t)  \, d \mathbf x 
\\[2ex]  \displaystyle 
+ \int_{\Omega}  \sum_{p=1}^n b_p(\mathbf x,t)
{\partial v \over \partial x_p}(\mathbf x,t) v(\mathbf x,t)
\, d \mathbf x  
\!+\! \int_{\Omega}  \hskip -1mm  (\lambda c \!-\! {1\over 2}  
c_t)(\mathbf x, t) |v(\mathbf x,t)|^2  \, d \mathbf x   = 0.  
\end{array}  
\end{eqnarray*}
Using uniform coercivity, the $L^\infty$ bounds,  and taking $\lambda$
large enough, we see that $\int_{\Omega} c(\mathbf x, t) |v(\mathbf x,t)|^2  
\leq \int_{\Omega} c(\mathbf x, 0) |v(\mathbf x,0)|^2   =0$. Therefore,
the solution is unique.

{\it Regularity.} Next, we differentiate with respect to $t$ to get
\begin{eqnarray} \begin{array}{l} \displaystyle
{d \over dt}
\int_{\Omega}  u_{t}(\mathbf x,t)  w(\mathbf x) \, c(\mathbf x,t) \, d \mathbf x 
+ \int_{\Omega}  \nabla u_t(\mathbf x,t)^T \mathbf b(\mathbf x,t) w(\mathbf x)
 \, d \mathbf x  \\ [2ex] \displaystyle
+ \int_{\Omega} \nabla u_t(\mathbf x,t)^T 
\mathbf A(\mathbf x,t) \nabla w(\mathbf x) \, d \mathbf x   
+ \int_{\Omega}  u_{t}(\mathbf x,t)  w(\mathbf x) \, c_t(\mathbf x,t) \, d \mathbf x    
\\ [2ex] \displaystyle
= \int_{\Omega}  f_t(\mathbf x,t)   w(\mathbf x)   \, d \mathbf x 
- \int_{\Omega}  u(\mathbf x,t)  w(\mathbf x) \, c_{tt}(\mathbf x,t) \, d \mathbf x
\\ [2ex] \displaystyle
+ \int_{\Omega}  \nabla u(\mathbf x,t)^T \mathbf b_t(\mathbf x,t) w(\mathbf x)
 \, d \mathbf x  + \int_{\Omega} \nabla u(\mathbf x,t)^T 
\mathbf A_t(\mathbf x,t) \nabla w(\mathbf x) \, d \mathbf x, 
\end{array} \label{vart_ref}
\end{eqnarray}
with $u_t(\mathbf x, 0) = w_0(\mathbf x)$. The functions
$\nabla u^T \mathbf b_t$, $\nabla u^T \mathbf A_t$, $u c_{tt}$
define linear forms in $H^1(\Omega)$.
Arguing as in Theorem 3.1, we see that the function $u_t$ is the unique
solution in $C([0,T];L^2(\Omega)) \cap L^2(0,T;H^1_{0}(\Omega))$  
of this problem. Then, (\ref{distributions}) implies that
$-{\rm div}(\mathbf A \nabla u) + \mathbf b^t \nabla u = -c u_t + f
\in C([0,T];L^2(\Omega))$  zero Dirichlet boundary condition. Elliptic 
regularity theory ensures that $u \in C([0,T];H^2(\Omega))$. 

{\it Stability.} The limiting solution inherits all the bounds established
on the approximating sequence. Therefore its $L^\infty([0,T];H^2(\Omega))$
and $H^1(0,T,H^1_0(\Omega))$ norms are bounded from above
in terms of constants depending on the parameters of the problem
and the norms of the data.
$\square$
\vskip 2mm

{\bf Theorem 3.2} {\it Under the hypotheses of Theorem 3.1, 
if $f\in L^q(\Omega \times [0,T])$ and $u_0 \in L^q(\Omega)$, then
$u$, its first and second order spatial derivatives,
and $u_t$ belong to $L^q(\Omega \times [0,T])$, $1<q<\infty$.}
\\
{\bf Proof.} 
We set $v/c = u$. Then $\nabla u = \nabla v/c - v/c^2 \nabla c$
and $cu_t = v_t - c_t/c v$. Therefore, $v$ is a solution of
\begin{eqnarray*}
v_t - {\rm div}\left({\mathbf A \over c} \nabla v \right) 
+  \mathbf A {\nabla c \over c^2} \nabla v + {\mathbf b^T \over c} 
\nabla v   
+ \left[ {\rm div}\left(\mathbf A {\nabla c \over c^2}\right)  
+  \mathbf b^T  {\nabla c \over c^2}   - {c_t \over c} \right] v = f.
\end{eqnarray*}
The result is a consequence of the regularity result stated in 
Theorem 9.1 in \cite{ladyzenskaya}.
$\square$ 
\vskip 2mm

\section{Well posedness results for the quasi-stationary submodels}
\label{sec:stsubmodels}

In this section we establish the pertinent existence and regularity results for
the elliptic submodels and the stationary transport problem in fixed domains.
Constructing solutions for the stationary transport problems considered
here is a non trivial issue. We are able to obtain them by a regularization
procedure under sign hypotheses on the velocity fields motivated by 
asymptotic studies, which will have to be preserved by any implemented 
scheme.

\subsection{Elliptic problems for displacements, velocities and concentrations}
\label{sec:elastic}

Consider the first the submodel for mechanical fields:
\begin{eqnarray} 
\begin{array}{ll} \displaystyle
\mu \Delta \mathbf u_s + (\mu +\lambda)  \nabla {\rm div}(\mathbf u_s) 
- \nabla  p =  \Pi \nabla \phi_s, & {\rm on } \; \Omega,  
\\[1ex] \displaystyle
\mu \Delta \mathbf v_s + (\mu +\lambda) \nabla {\rm div}(\mathbf v_s)  
= \nabla  p', 
& {\rm on } \; \Omega,  
\\[1ex] \displaystyle
k_h \Delta  p -{\rm div}(\mathbf v_s) =0,  
& {\rm on } \; \Omega,  
\\[1ex] \displaystyle
\Delta  p' =  (2\mu + \lambda) \Delta e', & {\rm on } \; \Omega,  
\\[1ex] \displaystyle
 p = p_{\rm ext}, \quad   p' = p_{\rm ext}'  & {\rm on }\; \Gamma, 
 \\[1ex] \displaystyle
\mathbf u  = 0, \quad \mathbf v  = 0,  & {\rm on } \; \Gamma_-, 
\\[1ex] \displaystyle
(\hat{\boldsymbol \sigma}(\mathbf u_s) - (p+\Pi \phi_s) \mathbf I) \mathbf n = \mathbf g, 
\quad  
(\hat{\boldsymbol \sigma}(\mathbf v_s) - p' \mathbf I) \mathbf n = \mathbf g',
& {\rm on } \; \Gamma_+. 
\end{array} \label{qselastic} 
\end{eqnarray}
We denote by $H^1_{0,-}(\Omega)$ the Sobolev space of $H^1(\Omega)$
functions vanishing on $\Gamma_-$. \vskip 1mm

{\bf Theorem 4.1.} {\it
Let $\Omega \subset \mathbb R^n$, $n=2,3,$ be an open bounded 
domain with  $C^4$ boundary $\partial \Omega$.  Let us assume
that $\phi_s \in H^1(\Omega)$ and $e' \in H^2(\Omega)$.
Given  positive constants $\mu$, $\lambda$, $k_h$, $\Pi$, there 
exists a unique solution 
$\mathbf u_s \in [H^2(\Omega)]^n \times [H^1_{0,-}(\Omega)]^n$,
$\mathbf v_s \in [H^3(\Omega)]^n \times [H^1_{0,-}(\Omega)]^n$,
$p \in H^4(\Omega)$,  $ p' \in H^2(\Omega)$  of (\ref{qselastic})
for any $p_{\rm ext}, p_{\rm ext}' \in \mathbb R$ and 
$\mathbf g, \mathbf g' \in \mathbb R^n$. \\
Moreover, if $\phi_s \in W^{1,q}(\Omega)$ and $e' \in W^{1,q}(\Omega)$,
$n<q<\infty$, then 
$p' \in W^{1,q}(\Omega)$, $\mathbf v_s \in W^{2,q}(\Omega)$,
$p \in W^{3,q}(\Omega)$ and $\mathbf u_s \in W^{2,q}(\Omega)$.
}

{\bf Proof.} The equation for $p'$ uncouples from the rest and provides a 
solution $p' \in H^2(\Omega)$ by classical theory for Laplace equations
\cite{brezis}. Next, the equation for $\mathbf v$ is a classical Navier elasticity
system which admits a unique solution $\mathbf v_s \in [H^2(\Omega)]^n 
\times [H^1_{0,-}(\Omega)]^n$ \cite{raviart}. Since the source $\nabla p' 
\in [H^1(\Omega)]^n$, elliptic regularity theory implies $\mathbf v_s \in [H^3(\Omega)]^n$. Now, 
${\rm div}(\mathbf v_s) \in H^2(\Omega)$ implies that the unique solution $p$ of 
the corresponding Poisson problem has $H^4(\Omega)$ regularity.
Finally, the equation for $\mathbf u_s$ is again a classical Navier elasticity
system with $L^2$ right hand side which admits a unique solution $\mathbf u_s \in [H^2(\Omega)]^n \cap [H^1_{0,-}(\Omega)]^n$.   \\
When $\phi_s \in W^{1,q}(\Omega)$ and  $e' \in W^{1,q}(\Omega)$,
we obtain the increased regularity \cite{kozlov}. Notice that since the boundary
values are constant, we can construct extensions to $H^k(\Omega)$ and
$W^{k,q}$ for the necessary $k, q$ \cite{brezis,raviart}. 
$\square$


\vskip 2mm
Now, the equation for the concentrations is:
\begin{eqnarray} \begin{array}{rcll}
-d  \Delta c  + {\rm div}  (\mathbf v_f c)  &=& - k_c g_c \phi_s, 
& \mathbf x \in \Omega, \\[0.5ex] \displaystyle
c &=& c_0     & \mathbf x \in \Gamma_-, \\[0.5ex] \displaystyle
{\partial c \over \partial \mathbf n} &=& 0 & \mathbf x \in \Gamma_+,
\end{array} \label{qschemical}  
\end{eqnarray}
given positive constants $d, c_0, k_c, g_c$ and known functions 
$\mathbf v_f$ and $\phi_s$.
\vskip 1mm

{\bf Theorem 4.2.} {\it Let $\Omega \subset \mathbb R^n$, $n=2,3,$ be an
open bounded domain with  $C^2$ boundary $\partial \Omega$.  Given 
positive constants $k_c, g_c, d, c_0$, a vector function 
$\mathbf v_l \in [H^1(\Omega)]^n \cap C(\overline{\Omega})$,
and a positive function $\phi_b \in  L^2(\Omega)$ there 
exists a unique nonnegative solution $c \in H^1(\Omega)$
of (\ref{qschemical})  provided $d$ is sufficiently large.  }

{\bf Proof.} Set $c= \tilde c + c_0$. The resulting problem admits the 
variational formulation: Find $\tilde c \in H^1_{0,-}(\Omega)$ such that
\begin{eqnarray*}
d \int_{\Omega}    \nabla \tilde  c^T  \nabla w \, d \mathbf x 
 -  \int_{\Omega} \mathbf v_f^T \tilde c \nabla w \, d \mathbf x 
+  \int_{\Gamma_+}   \tilde c w \mathbf v_l^T \mathbf n \, dS_{\mathbf x}  \\
 =  - k_c g_c \int_{\Omega} \phi_s w \, d \mathbf x + c_0 \int_\Omega
\mathbf v_f^T \nabla w  \, d \mathbf x,
\end{eqnarray*}
for all $ w \in H^1_{0,-}(\Omega)$. The continuous bilinear form is coercive 
provided $d$ is large enough compared to $\| \mathbf v_f\|_\infty$. Thus,
we have a unique solution $\tilde c \in H^1_{0,-}(\Omega)$ with 
$H^2(\Omega)$ regularity.

The function $c^- \in H^1_{0,-}(\Omega)$ satisfies
\begin{eqnarray*}
d \int_{\Omega}  |\nabla c^{-}|^2 \, d \mathbf x 
- \int_{\Omega} \mathbf v_f^T  c^{-} \nabla c^{-}  \, d \mathbf x 
+ \int_{\partial \Omega^+} |c^{-}|^2 \mathbf v_f^T
\mathbf n \, dS_\mathbf x   =
- k_c g_c \int_{\Omega}  \phi_s c^- \, d \mathbf x \leq 0.
\end{eqnarray*}
Coercivity  implies $c^{-}=0$ and $c \geq 0$  provided 
$d$ is large enough compared to $\| \mathbf v_l  \|_{\infty}$. 

For uniqueness, assume we have two positive solutions 
$c_1$ and $c_2$ in $H^1(\Omega)$ and set $c = c_1 - c_2
\in H^1_{0,-}(\Omega)$. 
Then $u$ is a solution of 
\begin{eqnarray*}
\begin{array}{rcll}
-  d \Delta c + {\rm div}  (\mathbf v_l c)   &=&  
0, & \mathbf x \in \Omega, \\[0.5ex] \displaystyle
c &=& 0, &  \mathbf x \in \partial \Omega^-, \\ [0.5ex] \displaystyle
{\partial c \over \partial \mathbf n} &=& 0, &  
\mathbf x \in \partial \Omega^+.  
\end{array}  
\end{eqnarray*}
The variational equation with test function $c$ and coercivity
imply $c=0$, that is, $c_1= c_2$. $\square$


\subsection{Conservation law for volume fractions}
\label{sec:conservation}

Consider the equation
\begin{eqnarray}
 {\rm div}(-\mathbf v_f  \phi_f) + k_s g_s \phi_f =  k_s g_s ,
 \quad   \mathbf x \in \Omega,   \label{vfractionl_stinfty}
\end{eqnarray}
where $k_s$ and  $g_s$ are positive constants and 
$\mathbf v_f$ a known function. \vskip 1mm

{\bf Theorem 4.3.} {\it
Let $\Omega \subset \mathbb R^n$, $n=2,3,$ be a thin open, bounded 
subset, with $C^4$ boundary $\partial \Omega$.  
Let $\mathbf v_f \in [H^2(\Omega) \cap C(\overline{\Omega})]^n$
such that  ${\rm div}(\mathbf v_f)\leq 0$ in $\Omega$,
${\rm div}(\mathbf v_f) \in L^\infty(\Omega) $
and $\mathbf v_f^T   \mathbf n \leq 0$ a.e. on $\partial 
\Omega$.  
We assume that $ \nabla \mathbf v_f \in [L^{\infty}(\Omega)]^{n^2}$ with 
$\|\nabla \mathbf v_f\|_{[L^{\infty}]^{n^2}}$ small enough compared to 
$k_s g_s$.
Then, given positive constants $k_s$ and  $g_s$,  there exists 
a  solution $\phi_f \in L^2(\Omega)$ of  (\ref{vfractionl_stinfty}) in the 
sense of distributions. Moreover,
\begin{itemize}
\item $0 \leq \phi_f \leq 1$ on $\Omega$ and $\phi$ does not vanish in 
sets of positive measure.
\item $\phi_f \in H^1(\Omega)$ is the unique solution of the
variational formulation in $H^1(\Omega)$ and 
\begin{eqnarray*}
{1\over 2} k_s g_s
 \| \nabla \phi  \|_{L^2} \leq 
\|\nabla {\rm div}(\mathbf v_f)\|_{[L^2]^n}.
\end{eqnarray*}
\item If we assume that $\Omega$  is a thin domain for which 
$\mathbf n \sim \mathbf e_n$ and ${\rm div}(\mathbf v_f) \in
W^{1,q}(\Omega)$, $n<q<\infty$,  then
$\nabla \phi_f \in L^q(\Omega)$ and
\begin{eqnarray*}
{1\over 2} k_s g_s
 \| \nabla \phi  \|_{L^q} \leq 
\|\nabla {\rm div}(\mathbf v_f)\|_{[L^q]^n}.
\end{eqnarray*}
\end{itemize}
}

{\bf Proof.}
{\it Existence.} For each $\varepsilon >0$, we follow \cite{beirao} and let 
$\phi_\varepsilon \in  H^1(\Omega)$ be  the solution of the variational
formulation
\begin{eqnarray*} b(\phi_\epsilon, w) =
\varepsilon \int_\Omega \nabla \phi _\varepsilon^T \nabla w \, d \mathbf x
+ \int_\Omega \mathbf v_f^T \phi _\varepsilon \nabla w \, d \mathbf x
- \int_{\partial \Omega} \phi_\varepsilon  w \, \mathbf v_f^T \mathbf n
d S_{\mathbf x} \\
+ \int_\Omega k_s g_s  \, \phi_\varepsilon w \, d \mathbf x  
= \int_\Omega k_s g_s  \, w \, d\mathbf x  = L(w), \quad \forall  w \in H^1(\Omega)
\end{eqnarray*}
of
\begin{eqnarray}
- \varepsilon \Delta  \phi _\varepsilon  - {\rm div}(\mathbf v_f  \phi_\varepsilon)
+ k_s g_s  \phi_\varepsilon = k_s g_s  \; {\rm in \,} \Omega, \quad 
{\partial  \phi_\varepsilon  \over \partial \mathbf n } = 0  \; {\rm on \,} 
\partial \Omega.
\label{eq_phie}
\end{eqnarray}
The bilinear form $b(\varphi, w)$ is continuous on $H^1(\Omega)$ 
\cite{raviart}, while the linear form $L$ is continuous on $L^2(\Omega).$ 
Since ${\rm div}(\mathbf v_f) \leq 0$ and $\mathbf v_f^T  \mathbf n \leq 0$, 
the bilinear form $b$ is also coercive in $H^1(\Omega)$. Indeed,
\begin{eqnarray*}
\int_\Omega \hskip -1mm  \mathbf v_f^T \phi_\varepsilon  \nabla 
\phi_\varepsilon  d \mathbf x
= {1\over 2}\int_\Omega \hskip -1mm  \mathbf v_f^T \nabla 
|\phi_\varepsilon|^2  d \mathbf x =  
{1\over 2}\int_{\partial \Omega} \hskip -1mm |\phi_\varepsilon|^2 
\mathbf v_f^T \mathbf n d \mathbf x
 - {1\over 2} \int_\Omega \hskip -1mm {\rm div}(\mathbf v_f) 
|\phi_\varepsilon|^2   d \mathbf x.
\end{eqnarray*}
The positive term $ - \int_\Omega {\rm div}(\mathbf v_f) |\phi_\varepsilon|^2 
d \mathbf x $ is finite because $|\phi_\varepsilon|^2 \in L^2(\Omega)$
thanks to Sobolev embeedings \cite{adams,brezis}. Since the bilinear form
$\varepsilon \int_\Omega \nabla \phi^T\nabla w \, d \mathbf x
+  \int_\Omega k_s g_s \phi w \, d \mathbf x$
is coercive in $H^1(\Omega)$, we have a  unique solution 
$\phi_\varepsilon \in H^1(\Omega)$ by Lax Milgram's theorem \cite{brezis}.
We set $w=\phi _\varepsilon$ and apply  Young's inequality \cite{brezis} 
to obtain the uniform bound
$  \| \phi_\varepsilon \|_{L^2} \leq {\rm meas}(\Omega)^{1/2}$ from
\begin{eqnarray*} 
0 \leq \varepsilon \int_\Omega |\nabla \phi _\varepsilon|^2 \, d \mathbf x
- {1\over 2} \int_{\partial \Omega}  |\phi_\varepsilon|^2  \mathbf v_f^T 
\mathbf n \, d S_{\mathbf x} +  \int_\Omega \left[ - {1\over 2} 
{\rm div}(\mathbf v_f)  + k_s g_s \right] |\phi_\varepsilon|^2  \, d \mathbf x  \\ 
= \int_\Omega k_s g_s \phi_\varepsilon  \, d\mathbf x   
\leq \| k_s g_s \|_{L^2} \left( \int_\Omega  |\phi_\varepsilon|^2 \right)^{1/2}.
\end{eqnarray*}
Each of the positive terms in the left hand side of the above 
inequality are uniformly bounded too. 
Thus, we can extract a subsequence $\phi_{\varepsilon'}$ such that 
$ \phi_{\varepsilon'}$ tends weakly in $L^2(\Omega)$ to a limit $  \phi$,
and $\varepsilon \nabla \phi_\varepsilon$ tends strongly to zero.
Setting $w \in C_c^\infty(\Omega)$ in the variational formulation,
and taking limits \cite{amm16,lions}, $\phi$ is a solution of (\ref{vfractionl_stinfty}) 
in the sense of distributions. The variational equation holds with $\epsilon =0$, 
replacing the boundary integral by the duality  
$_{H^{-1/2}(\partial \Omega)}<\phi \, \mathbf v_f^T 
\mathbf n, w>_{H^{1/2}(\partial \Omega)}$ for $w\in H^1(\Omega)$
\cite{bamberger}.

{\it $L^\infty$ estimates. }
Setting $\psi_\varepsilon =  \phi_\varepsilon - 1$ and  $w = \psi_\varepsilon^+$
we get
\begin{eqnarray*} 
\varepsilon \int_\Omega |\nabla \psi _\varepsilon^+|^2 \, d \mathbf x - 
{1\over 2} \int_{\partial \Omega}  |\psi_\varepsilon^+|^2 \, 
\mathbf v_f^T \mathbf nd S_{\mathbf x} 
+  \int_\Omega \left[ - {1\over 2} {\rm div}(\mathbf v_f)  +
k_s g_s \right] |\psi_\varepsilon^+|^2  \, d \mathbf x    \\ 
= \int_\Omega  {\rm div}(\mathbf v_f) \psi_\varepsilon^+ d\mathbf x \leq 0.
\end{eqnarray*}
Thus, $\psi_\varepsilon^+=0$ and $ \phi_\varepsilon \leq 1$. Similarly,
we set  $\psi_\varepsilon =  - \phi_\varepsilon$  to find
\begin{eqnarray*} 
\varepsilon \int_\Omega |\nabla \psi _\varepsilon^+|^2 \, d \mathbf x - 
{1\over 2} \int_{\partial \Omega} (\mathbf v_f^T \mathbf n) |\psi_\varepsilon^+|^2 
\, d S_{\mathbf x} +  \int_\Omega \left[ - {1\over 2} {\rm div}(\mathbf v_f)  +
k_s g_s \right] |\psi_\varepsilon^+|^2  \, d \mathbf x    \\ 
= - \int_\Omega  k_s g_s  \psi_\varepsilon^+ d\mathbf x \leq 0.
\end{eqnarray*}
Thus, $\psi_\varepsilon^+=0$ and $ \phi_\varepsilon \geq 0.$
Weak limits $ \phi$ in $L^2$ inherit these  properties.
Moreover, (\ref{vfractionl_stinfty}) implies that $ \phi$ cannot vanish in
sets of positive measure.

{\it $H^1$ Regularity.}
Elliptic regularity for system (\ref{eq_phie}) implies that  
$\phi_\varepsilon \in H^2(\Omega)$ \cite{adn2,brezis}. We multiply 
(\ref{eq_phie}) by $\Delta \phi_\varepsilon$ and  integrate over $\Omega$ 
to get
\begin{eqnarray*}
-  \varepsilon \int_\Omega |\Delta \phi _\varepsilon|^2  d \mathbf x
- \int_\Omega \!\! \mathbf v_b^T \nabla \phi_\varepsilon \Delta 
\phi _\varepsilon d \mathbf x
+ \int_\Omega \!\! \left[ - {\rm div}(\mathbf v_b) + k_s g_s
\right] \phi_\varepsilon  \Delta \phi_\varepsilon  d \mathbf x  \\
= \int_\Omega k_s g_s \Delta \phi _\varepsilon d \mathbf x.
\end{eqnarray*}
Integrating by parts, and using the boundary condition, we find
\begin{eqnarray*}
-  \varepsilon \int_\Omega |\Delta \phi _\varepsilon|^2  d \mathbf x
+ \int_\Omega \left[  {1\over 2} {\rm div}(\mathbf v_f)  -
k_s g_s \right] |\nabla \phi_\varepsilon|^2  d \mathbf x
+ {1\over 2} \int_{\partial \Omega}  |\nabla \phi_{\varepsilon}|^2 \mathbf v_f^T 
\mathbf n d S_{\mathbf x} = \\
\int_\Omega \nabla \! \left[ - {\rm div}(\mathbf v_f)   + k_s g_s \right] ^T
\phi_\varepsilon  \nabla \phi_\varepsilon  d \mathbf x  
 -   \int_\Omega v_{l,j,x_k} \phi_{\varepsilon, x_j}  \phi_{\varepsilon, x_k}  
 d \mathbf x. 
\end{eqnarray*}
We know that $0\leq \phi _\varepsilon \leq 1$.
Therefore, 
\begin{eqnarray*}
\int_\Omega \left[  -{1\over 2} {\rm div}(\mathbf v_f)  +
k_s g_s  \right] |\nabla \phi_\varepsilon|^2  d \mathbf x
\leq  \|\nabla {\rm div}(\mathbf v_f)\|_{[L^2]^n} 
 \| \nabla  \phi_\varepsilon \|_{L^2} \\
+  \int_\Omega |v_{l,j,x_k} \phi_{\varepsilon, x_j}  \phi_{\varepsilon, x_k}| 
d \mathbf x.
\end{eqnarray*}
If $\|\nabla \mathbf v_{l} \|_{[L^\infty]^{n^2}}$ is small enough compared
to $k_s g_s $
\begin{eqnarray*}
{1\over 2} k_s g_s   \| \nabla \phi_\varepsilon \|_{L^2} \leq 
\|\nabla {\rm div}(\mathbf v_f)\|_{[L^2]^n}.
\end{eqnarray*}
We extract a subsequence $\phi_{\varepsilon'}$ converging weakly in 
$H^1(\Omega)$ to a limit $\phi$, strongly in $L^2(\Omega)$, and 
pointwise in $\Omega$. The traces of $\phi$  on $\partial \Omega$ 
belong to $L^2(\partial \Omega)$, and are weak limits of traces
of  $\phi_{\varepsilon'}$. Passing to the limit in the variational formulation
for (\ref{eq_phie}), $\phi \in H^1(\Omega)$ is a solution with 
$\epsilon =0$ which inherits these bounds.

{\it Uniqueness.} Given two solutions $\phi_1, \phi_2  \in H^1(\Omega)$, 
we set $\psi = \phi_1-\phi_2$. Subtracting the variational equations 
we get for the test function $\psi \in H^1(\Omega)$  
\begin{eqnarray*} 
- {1\over 2} \int_{\partial \Omega} (\mathbf v_f^T \mathbf n) |\psi|^2 
\, d S_{\mathbf x} +  \int_\Omega \left[ - {1\over 2} {\rm div}(\mathbf v_f)  +
k_s g_s  \right] |\psi|^2  \, d \mathbf x  = 0,
\end{eqnarray*}
that is, $\phi_1=\phi_2$ in view of the signs. $\square$ 

{\it $W^{1,q}$ regularity.} 
By elliptic regularity,  $\phi _\varepsilon \in W^{3,q}(\Omega)$, 
since the source in (\ref{eq_phie}) belongs to $W^{1,q}(\Omega)$. 
Following \cite{beirao}, we differentiate (\ref{eq_phie}) with respect 
to $x_k$, multiply by $h(\phi_\varepsilon) \phi_{x_k}$ for 
$h(\phi_\varepsilon) = (|\nabla \phi_\varepsilon|^2 + \delta)^{(q-2)/2}$,
add $k$ and integrate over $\Omega$ to get
\begin{eqnarray*}
- \varepsilon \int_{\Omega} \Delta(\nabla \phi_\varepsilon)^T 
h(\phi_\varepsilon) \nabla \phi_\varepsilon \, d \mathbf x + 
\int_{\Omega} k_s g_s h(\phi_\varepsilon) |\nabla \phi_\varepsilon|^2 \, 
d \mathbf x  \\
 -  \int_{\Omega} \!\! v_{l,i} \phi_{\varepsilon,x_i x_k} h(\phi_\varepsilon) 
\phi_{\varepsilon, x_k}   d \mathbf x   
- \int_{\Omega}    v_{l,i,x_k} \phi_{\varepsilon, x_i} h(\phi_\varepsilon)
\phi_{\varepsilon, x_k}     d \mathbf x   \\
-  \int_{\Omega}   {\rm div}(\mathbf v_f) h(\phi_\varepsilon) |\nabla 
\phi_\varepsilon|^2  d \mathbf x   -
\int_{\Omega} \nabla({\rm div}(\mathbf v_f))^T h(\phi_\varepsilon) 
\phi_\varepsilon \nabla \phi_\varepsilon \, d \mathbf x = 0.
\end{eqnarray*}
Sum over repeated indices is intended. Notice that Lemma 3.1 from \cite{beirao}
holds in our framework for our thin domains, so that the first term is nonnegative.
The fourth term becomes
\begin{eqnarray*}
{1\over q} \int_{\Omega} {\rm div}(\mathbf v_f)(|\nabla \phi_\varepsilon|^2 
+ \delta)^{q/2} d \mathbf x 
- {1\over q} \int_{\partial \Omega} (|\nabla \phi_\varepsilon|^2 
+ \delta)^{q/2}  \mathbf v_{l}^T \mathbf n \, dS_{\mathbf x}.
\end{eqnarray*}
Putting all together we get
\begin{eqnarray*}
 \int_{\Omega} k_s g_s h(\phi_\varepsilon) |\nabla \phi_\varepsilon|^2 \, 
d \mathbf x 
\leq - {1\over q} \int_{\Omega} {\rm div}(\mathbf v_f)(|\nabla \phi_\varepsilon|^2 
+ \delta)^{q/2} d \mathbf x \\
+ \int_{\Omega}  v_{l,i,x_k} \phi_{\varepsilon, x_i} \phi_{\varepsilon, x_k}  
h(\phi_\varepsilon) \, d \mathbf x
+ \int_{\Omega} {\rm div}(\mathbf v_f) h(\phi_\varepsilon) |\nabla \phi_\varepsilon|^2 
\, d \mathbf x \\
+ \int_{\Omega} \nabla({\rm div}(\mathbf v_f))^T h(\phi_\varepsilon) 
\phi_\varepsilon \nabla \phi_\varepsilon \, d \mathbf x.
\end{eqnarray*}
We let $\delta \rightarrow 0$ and use that
$\| \nabla \mathbf v_f \|_{[L^\infty]^{n^2}}$ is small enough to find
\begin{eqnarray*}
{1 \over 2} k_s g_s  \int_{\Omega}  |\nabla \phi_\varepsilon|^q \, d \mathbf x 
\leq  \| \nabla({\rm div}(\mathbf v_f)) \|_{L^q}    \| |\nabla \phi_\varepsilon| \|_{L^{q}}^{q-1},
\end{eqnarray*}
which yields the bound we seek letting $\varepsilon \rightarrow 0$. 
$\square$ 


\section{Well posedness results for the full model with a known boundary 
dynamics}
\label{sec:full}

Once we have analyzed the different submodels, we consider the whole 
system when the boundary of the domains $\Omega^t$ moves with time
according to a given dynamics
\begin{eqnarray} 
\begin{array}{ll} \displaystyle
\mu \Delta \mathbf u_s + (\mu +\lambda)  \nabla {\rm div}(\mathbf u_s) 
- \nabla  p =  \Pi \nabla \phi_s, 
& {\rm in } \; \Omega^t,  
\\[1ex] \displaystyle
\mu \Delta \mathbf v_s + (\mu +\lambda) \nabla {\rm div}(\mathbf v_s)  
= \nabla  p' , 
& {\rm in } \; \Omega^t,  
\\[1ex] \displaystyle
k_h \Delta  p = {\rm div}(\mathbf v_s),  
& {\rm in } \; \Omega^t,  
\\[1ex] \displaystyle
 \Delta  p' =  (2\mu + \lambda) \Delta e', & {\rm in } \; \Omega^t,  
\\[1ex] \displaystyle
 p = p_{\rm ext}, \quad   p = p_{\rm ext}'  & {\rm on }\; \Gamma^t, 
 \\[1ex] \displaystyle
\mathbf u_s = 0, \quad \mathbf v_s = 0,  & {\rm on } \; \Gamma_-^t, 
\\[1ex] \displaystyle
(\hat{\boldsymbol \sigma} (\mathbf u_s) - (p+\Pi \phi_s) \mathbf I) 
\mathbf n = \mathbf g, & {\rm on } \; \Gamma_+^t, 
\\[1ex] \displaystyle
(\hat{\boldsymbol \sigma} (\mathbf v_s) - p' \mathbf I) 
\mathbf n  = \mathbf g'(\nabla \mathbf u_s), & {\rm on } \; \Gamma_+^t, 
\end{array} \label{felastic} \\[2ex] 
\begin{array}{ll} \displaystyle
{\rm div}  (-{\mathbf v}_f \phi_f) + k_s g_s \phi_f = k_s g_s, 
\hskip 1cm &  {\rm in } \; \Omega^t,  
\\[1ex] \displaystyle
\mathbf v_f  = - \xi_\infty \nabla p +  \mathbf v_s, 
\quad \phi_f+\phi_s = 1, \hskip 1cm &  {\rm in } \; \Omega^t,
\end{array} \label{ffraction} \\[2ex] 
\begin{array}{ll} \displaystyle
{de' \over dt} = k_h (2 \mu + \lambda)  \Delta e', 
\hskip 2.9cm & \; \mbox{in $\Omega^t$},  
\\[1ex] \displaystyle
e' = e_{\rm ext},  & \hskip 1mm {\rm on }\; \Gamma^t,
\\[1ex] \displaystyle
e'(0) = e_0,  & \hskip 1mm {\rm on }\; \Omega^0,
\end{array} \label{fdiv} \\[2ex] 
\begin{array}{ll}
-d  \Delta c  + {\rm div}  (\mathbf v_f c)  = - k_c g_c \phi_s,  \hskip 1.8cm
&  {\rm in } \; \Omega^t, \\[0.5ex] \displaystyle
c = c_0     &    {\rm on } \; \Gamma^t_-, \\[0.5ex] \displaystyle
{\partial c \over \partial \mathbf n} = 0 &     {\rm on } \; \Gamma^t_+.
\end{array} \label{fchemical}    
\end{eqnarray}

\vskip 1mm
{\bf Theorem 5.1.} {\it Let $\Omega^t \subset \Omega \subset \mathbb R^n$, 
$n=2,3$, $t \in [0,T]$, be a family of open bounded $C^4$ domains. The lower  boundary $\Gamma_-$ is fixed, while the upper boundary $\Gamma_+^t$ is obtained deforming $\Gamma^0_+$ along a vector field $\mathbf \nu(\mathbf x) 
\in C(\overline{\Omega}) \cap C^4(\Omega)$. Assume that 
\begin{itemize}
\item $e_{\rm ext}(t)$,  $\mathbf g(t)$, $\mathbf g'(t)$, $p_{\rm ext}(t)$, 
$p_{\rm ext}'(t)$, $c_0(t) \in C([0,T])$, $e_0 \in L^2(\Omega^0) \cap
L^q(\Omega^0)$, for $q>n$,
\item $e_{\rm ext}$, $\mathbf g'$, $\Pi$ and $p_{\rm ext}$ are small enough.
\end{itemize}
Given positive constants $\mu$, $\lambda$, $\Pi$, $k_h$, $k_s$, 
$k_c$, $g_s$, $g_c$, $\xi_\infty$, and $d$ large enough,
system (\ref{felastic})-(\ref{fchemical}) admits a unique solution 
$e'  \in  H^2(\Omega^t)\cap W^{2,q}(\Omega^t)$, 
$\mathbf u_s  \in [H^2(\Omega^t)]^n \cap W^{2,q}(\Omega^t)$,
$\mathbf v_s, \mathbf v_f  \in [H^3(\Omega^t)]^n$,
$p  \in  H^4(\Omega^t)$, $p'  \in  H^2(\Omega^t)$,
$\phi_f, \phi_s \in H^1(\Omega^t)\cap W^{1,q}(\Omega^t)$,
$c \in H^2(\Omega^t)$, for $q >n$,
satisfying $c \geq 0$ and $0 \leq \phi_f, \phi_s \leq 1$, $t \in [0,T]$.
Moreover, the norms of the solutions are bounded in terms
of the parameters and data of the problem.}

{\bf Proof.} 
Assume first that $\mathbf g'(\nabla \mathbf u_s)$ does not depend
on $\mathbf u_s$. Then, the result is a consequence of Corollary 3.3, 
Theorems 4.1-4.3 and Sobolev embeddings 
\cite{adams} (neither $L^q$ regularity nor conditions on the domain
geometry nor smallness assumptions are needed). 
We calculate the unknowns according to the 
sequence $e'$, $p'$, $\mathbf v_s$, $p$, $\mathbf v_f$, $\phi_f$, 
$\phi_s$, $\mathbf u_s$, and $c$.  

When $\mathbf g'(\nabla \mathbf u_s)$ does depend on 
$\mathbf u_s$, we construct $e'$ thanks to Corollary 3.3. For 
each fixed $t>0$, $e' \in H^2(\Omega^t)\cap W^{2,q}(\Omega^t)$
and we can construct $p' \in H^2(\Omega^t)\cap W^{2,q}(\Omega^t)$. 
Next, we solve the quasi-stationary system by means of an 
iterative scheme.
At each step $\ell$, we freeze $\Pi \nabla \phi_s^{(\ell -1)}$ in 
the equation for $\mathbf u_s^{(\ell)}$ and 
$\mathbf g'(\nabla \mathbf u_s^{\ell-1})$ in the boundary
condition for $\mathbf v_s^{(\ell)}$. Initially, 
we set $\phi_s^{(0)}= \phi_\infty \in (0,1)$ constant and 
$\phi_f^{(0)}=1-\phi_\infty$. We set $\mathbf u^{(0)}=0$.
Theorem 4.1, Theorem 4.2, Theorem 4.3 guarantee
the existence of  $\mathbf v_s^{(1)}$, $p^{(1)}$, 
$\mathbf u_s^{(1)}$, $\mathbf v_f^{(1)}$, $\phi_f^{(1)}$, 
$\phi_s^{(1)}$, and $c^{(1)}$, with the stated regularity. 

In a similar way, given all the fields at step $\ell-1$, we can 
construct the solutions for step $\ell$. 
Notice that $\mathbf v_f^{(\ell-1)} \in W^{2,q} $ implies
$\mathbf v_f^{(\ell-1)}  \in W^{1,\infty}(\Omega)$ and 
$\mathbf v_f^{(\ell-1)}  \in C(\overline{\Omega})$. To apply Theorem
4.3 we also need to satisfy smallness and sign assumptions
that we will consider later. Assuming they hold, we get for the elliptic 
system  involving $\mathbf v_s^{(\ell)}$, $\mathbf u_s^{(\ell)}$, 
$p^{(\ell)}$ and for the transport equation for $\phi_s^{(\ell)}$
\begin{eqnarray*}
\| p^{(\ell)} \|_{H^2(\Omega^t)} 
+ \| \mathbf v_s^{(\ell)} \|_{H^2(\Omega^t)} +
\| \mathbf u_s^{(\ell)} \|_{H^2(\Omega^t)} \leq          
C_1^t [\Pi \| \nabla \phi_s^{(\ell-1)} \|_{L^2(\Omega^t)}             \\ 
+ \| \nabla p' \|_{L^2(\Omega^t)} + \|p_{\rm ext} \|_{H^{3/2}(\Gamma^t_+)}
+ \|\mathbf g'(\nabla \mathbf u_s^{(\ell-1)}) \|_{H^{1/2}(\Gamma^t_+)}
+ \|\mathbf g \|_{H^{1/2}(\Gamma^t_+)} ],                          \\[1ex]
\| p^{(\ell)} \|_{W^{2,q}(\Omega^t)} 
+ \| \mathbf v_s^{(\ell)} \|_{W^{2,q}(\Omega^t)} +
\| \mathbf u_s^{(\ell)} \|_{W^{2,q}(\Omega^t)} \leq    
C_2^t [ \Pi \| \nabla \phi_s^{(\ell-1)} \|_{L^q(\Omega^t)} +       \\     
\| \nabla p' \|_{L^q(\Omega^t)} +  
\| p_{\rm ext} \|_{W^{1-{1\over q},q}(\Gamma^t_+)}
+ \|\mathbf g'(\nabla \mathbf u_s^{(\ell-1)}) \|_{W^{1-{1\over q},q}(\Gamma^t_+)}
+ \|\mathbf g \|_{W^{1-{1\over q},q}(\Gamma^t_+)} ],         \\[0.5ex]
\| p^{(\ell)} \|_{W^{3,q}(\Omega^t)} \leq C^t_3 [
\| \mathbf v_s^{(\ell)} \|_{W^{1,q}(\Omega^t)} +
\| p_{\rm ext} \|_{W^{3-1/q,q}(\Gamma^t)} ]                        \\[1ex]
\| \mathbf v_f^{(\ell)} \|_{W^{2,q}(\Omega^t)} \leq 
\xi_\infty \|  p^{(\ell)} \|_{W^{3,q}(\Omega^t)} +
\| \mathbf v_s^{(\ell)} \|_{W^{2,q}(\Omega^t)}                     \\[1ex]
{1\over 2} k_s g_s \| \nabla \phi_f^{(\ell)}  \|_{L^q} \leq 
\|\nabla {\rm div}(\mathbf v_f^{(\ell)})\|_{[L^q]^n}.                                
\end{eqnarray*}
Notice that $\nabla \phi_f^{(\ell)} = - \nabla \phi_s^{(\ell)}$.
Combining the above inequalities, and provided $\Pi$ and 
$\mathbf g'$ are small enough, we obtain an upper bound for
$\| \mathbf v_f^{(\ell)} \|_{W^{2,q}(\Omega^t)} $,
$\| \mathbf v_s^{(\ell)} \|_{W^{2,q}(\Omega^t)} $,
$\| p^{(\ell)} \|_{W^{2,q}(\Omega^t)} $,
$\| \phi_s \|_{W^{1,q}(\Omega^t)} $, in terms of 
constants depending on the problem data and parameters,
and also on time, but remain bounded in time for $t\in [0,T]$.



We guarantee by induction the smallness of 
$\|\mathbf v_f^{(\ell)}|_{[W^{1,\infty}]}$ and 
${\rm div}(\mathbf v_f^{(\ell)}) \leq 0$, $\mathbf v_f^{(\ell)} \cdot 
\mathbf n \leq 0$.
Initially, $\phi_s^{(0)}$ is constant and $\nabla \phi_s^{(0)}=0$. We 
construct $\mathbf v_s^{(1)}$ and $p^{(1)}$ in such a way that 
$\|\mathbf v_s^{(1)} \|_{[W^{2,q}]^n}$,
$\|p^{(1)} \|_{[W^{3,q}]^n}$ and  $\|\mathbf v_f^{(1)} \|_{[W^{2,q}]^n}$ 
are bounded in terms of the problem parameters and data. By Sobolev
injections for $n < q < \infty$, $\|\mathbf v_s^{(1)} \|_{[W^{1,\infty}]^n}$ 
satisfies  a similar bound, and can be made as small as required by 
making $\mathbf g'$ and $p_{\rm ext}$ small. Then,
$\|\nabla \phi_f^{(1)}\|_{L^q}$   is bounded by 
$\|\mathbf v_f^{(1)} \|_{[W^{2,q}]^n}$  and is equally small.
Furthermore, ${\rm div}(\mathbf v_f^{(1)}) \phi_f^{(1)} + \mathbf v_f^{(1)} 
\nabla \phi_f^{(1)} = - k_s g_s \phi_f^{(1)} \leq 0$. Since
$\mathbf v_f^{(1)} $ and $\nabla \phi_f^{(1)} $ are small compared
to $- k_s g_s \phi_f^{(1)} \leq 0$ which is
almost constant. Thus, ${\rm div}(\mathbf v_l^{(1)}) \leq 0$.
Finally, $\int_{A} {\rm div}(\mathbf v_l^{(1)}) \, d {\mathbf x}
= \int_{\partial A} \mathbf v_l^{(1)} \cdot \mathbf n \, d {S_{\mathbf x}}
\leq 0$ for all $A \subset \Omega$ so that
$\mathbf v_l^{(1)} \cdot \mathbf n \leq 0$ on $\partial \Omega$.

By induction, if $\|\mathbf v_f^{(\ell-1)} \|_{[W^{1,\infty}]^n}$ is small
and $\mathbf v_f^{(\ell-1)}$ satisfies the sign conditions, we 
can repeat the argument to show that this holds for $\mathbf v_f^{(\ell)}$
too and that it also satisfies the sign conditions. We need to estimate 
$\| \nabla {\rm div}(\mathbf v_f^{(\ell-1)}) \|_{[L^q]^n}$, which is possible
since $\Pi$ is small.

These estimates allow us to extract subsequences converging weakly
to limits $\mathbf v_s$, $\mathbf u_s$, $p$,  $\phi_s$ satisfying
variational formulations of the equations. Problem (\ref{fchemical})
is already studied in Theorem 4.2.
$\square$

A similar result (except for the uniqueness) can be obtained by means
of an iterative scheme if we allow for almost constant  
smooth coefficients $k_h(\phi_f)$ $\xi_\infty(\phi_f)$, $g_s(c)$, $g_s(c)$. 
\vskip 1mm

\section{Discussion and conclusions}
\label{sec:conclusions}
 
The study of biological aggregates and tissues often leads to complex
mixture models, combining transport equations for volume fractions of
different phases, with continuum models for mechanical behavior of
the mixture and chemical species \cite{kapellos,fstissue,fsbrain}. 
These models are set in domains that change with time, because cells 
grow, die and move and because of fluid transport within the biological
network.  Here, we have considered a fluid-solid mixture 
description of the spread of cellular systems called biofilms, which could 
be adapted to general tissues. These models involve different time scales, 
so that part of the equations are considered quasi-stationary, that is, they 
are stationary  problems solved at different times in different domains and 
with some time dependent coefficients. Such equations are coupled to time dependent problems set in moving domains and to variables not directly characterized by means of equations.

In this paper, we have developed mathematical frameworks to tackle
some of the difficulties involved in the construction of solutions for these
multiphysics systems and the study of their behavior.
First, we have  shown how to improve these models by characterizing time 
derivatives of solutions of stationary boundary value problems with varying coefficients set in moving domains in terms of complementary boundary 
value problems derived for them.  In this way we obtain a quasi-stationary
elliptic system for the mechanical variables of the solid phase, not only displacements and pressure, but also velocity, that can be solved at each time  coupled to the other submodels. This option is more stable 
than evaluating velocities as quotients of differences of displacements 
calculated in meshes of different spatial domains. On one side, the error committed is easier to control. On the other side, the computational is
cost smaller, since we use a single mesh at each time.
Once we know the velocity of the solid phase and the pressure, the velocity 
of the fluid phase follows by a Darcy type law. 
Next, we have devised an strategy to construct solutions of an auxiliary 
class of time  dependent linear diffusion problems set in moving domains 
with parametrizations satisfying a number of conditions. We are able to 
refer the model to a fixed domain and then solve by Galerkin type schemes.
The complete model involves a quasi-stationary transport problem. We
show that we can construct smooth enough solutions  by a regularization
procedure, under sign hypothesis on the fluid velocity field suggested by
asymptotic solutions constructed in simple geometries.
Once we know how to construct stable solutions of each submodel
satisfying adequate regularity properties, an iterative scheme allows us
 to solve the full problem when the time evolution of the boundary of the 
 spatial region occupied by the biological film is known.

In applications one must couple these models with additional lubrication
type equations for the motion of the film boundary, see equation (\ref{iheight}).
Perturbation analyses \cite{entropy} provide approximate
solutions with selfsimilar dynamics for $h$. Establishing existence and regularity 
results for such complex models that can guide construction of reliable numerical solutions is a completely open problem. 
The techniques we have developed are general and can be
applied in models with a similar  structure arising in other biological 
and chemical engineering applications.




\section*{Appendix: The model equations}
\label{sec:model}

We study biofilms as solid-fluid mixtures, composed of a solid biomass 
phase and a liquid phase formed by water carrying dissolved chemicals 
(nutrients, autoinducers, waste).
Under the equipresence hypothesis of mixtures, each location 
$\mathbf x$  in a biofilm can contain both phase simultaneously, 
assuming that no voids or air bubbles form inside. Let us denote by
$\phi_s(\mathbf x,t) $ the volume fraction of solid and by
$\phi_f(\mathbf x,t) $ the volume fraction of fluid, which satisfy
\begin{eqnarray}
\phi_s + \phi_f = 1. \label{fraction}
\end{eqnarray}
Taking the densities  the mixture and both constituents to be 
constant and equal to that of water $\rho_f= \rho_s= \rho= \rho_w$, 
the mass balance laws for  $\phi_s$ and $\phi_f$ are  
\cite{lanir, seminara}
\begin{eqnarray} 
{\partial \phi_s \over \partial t} + {\rm div}  (\phi_s {\mathbf v}_s) = 
 r_s(\phi_s,c),  \quad r_s(\phi_s,c) =
 k_s {c \over c + K_c} \phi_s, \label{lcsolid}  \\
 {\partial \phi_f \over \partial t} + {\rm div} (\phi_f {\mathbf v}_f) 
 = - r_s(\phi_s,c),   \label{lcfluid} 
\end{eqnarray}
where $\mathbf v_s$ and $\mathbf v_f$ denote the velocities of 
the solid and fluid components, respectively, 
$c$ is the substrate concentration and $r_s(\phi_s,c) =
 k_s {c \over c + K_c} \phi_s $ stands for the production of 
biomass due to nutrient consumption. The parameters $K_c$ 
(starvation threshold) and $k_s$ (intake rate) are positive 
constants.
The substrate concentration $c$ \cite{entropy,ibm} is governed by:
\begin{eqnarray}
{\partial  c \over \partial t} + {\rm div}  (\mathbf v_f c)  
- {\rm div} (d \nabla c)   =  -r_n(\phi_s,c), \quad
r_n(\phi_s, c) \!=\! \phi_s k_c { c \over c + K_c},
\label{lcnutrient} 
\end{eqnarray}
where  $r_n(\phi_s,c)$ represents consumption by the biofilm.
The parameters $d$ (diffusivity), $k_c$ (uptake rate) and 
$K_c$ (half-saturation) are positive constants. We impose
zero-flux boundary conditions on the air--biofilm interface
and constant  Dirichlet boundary condition on the agar--biofilm
interface.
In equation 
(\ref{lcnutrient}), typical
parameter values are such that the time derivatives
can be neglected. The solutions depend on time though the motion
of the biofilm boundary.

Adding up equations (\ref{lcsolid}) and (\ref{lcfluid}), we obtain
a conservation law for the growing mixture:
\begin{eqnarray}
0= {\rm div}  (\phi_s {\mathbf v}_s+ \phi_f {\mathbf v}_f) = 
{\rm div}  (\mathbf v) = {\rm div} (\mathbf v_s + \mathbf q),
\label{balancemixture}
\end{eqnarray}
where $\mathbf v= \phi_s {\mathbf v}_s+ \phi_f {\mathbf v}_f$ is the
averaged velocity  and
\begin{eqnarray}
\mathbf q = \phi_f (\mathbf v_f - \mathbf v_s) 
\label{flux}
\end{eqnarray}
is the filtration flux.

The theory of mixtures hypothesizes that the motion of each phase
obeys the usual momentum balance equations \cite{lanir}. In the absence 
of external body forces, the momentum balance for the solid and the
fluid reads
\begin{eqnarray} 
\rho \phi_s a_s + {\rm div} \boldsymbol \sigma_s 
+ \rho \phi_s (\mathbf  f_s + \nabla \pi_s) = 0,\quad 
\rho \phi_f a_f + {\rm div} \boldsymbol \sigma_f  
+ \rho \phi_f  (\mathbf  f_f + \nabla \pi) = 0. \label{motion}
\end{eqnarray}
In biofilms, the velocities $\mathbf v_s$ and $\mathbf v_f$  are small enough 
for inertial forces to be neglected, that is,
$\rho_s \mathbf a_s \approx \rho_f \mathbf a_f \approx \rho \mathbf a \approx 0$, 
where $\mathbf a_s, \mathbf a_f, \mathbf a$ denote the solid, fluid, and 
average accelerations.

Let us detail now expressions for the stresses and forces appearing
in these equations, following \cite{entropy,lanir}.
When the biofilm contains a large number of small pores, the 
stresses in the fluid are  
\begin{eqnarray}
\boldsymbol \sigma_f =  - \phi_f \, p \mathbf I,
\label{stressfluid}
\end{eqnarray}
$p$ being the pore hydrostatic pressure. In case large regions filled 
with fluid were present, the standard stress law for viscous fluids 
should be considered.
Under small deformations, and assuming an isotropic solid, the stresses 
in  the solid biomass are
\begin{eqnarray}
\boldsymbol \sigma_s =  \hat{\boldsymbol \sigma}_s  - \phi_s \, p \mathbf I, 
\quad \hat{\boldsymbol \sigma}_s = \lambda {\rm Tr} 
(\boldsymbol  \varepsilon({\mathbf u}_s))\, \mathbf I +  
2 \mu \, \boldsymbol  \varepsilon({\mathbf u}_s), \quad
\varepsilon_{ij}({\mathbf u})= {1\over 2} \Big(  {\partial u_i \over \partial x_j } 
+ {\partial u_j \over \partial x_i}  \Big),
\label{stresssolid}
\end{eqnarray}
where ${\mathbf u}_s$  is the displacement vector of the solid,
$ \boldsymbol  \varepsilon({\mathbf u})$ the deformation tensor,
and $\lambda, \mu,$ the Lam\'e constants.
The stresses in the solid are due to interaction with the fluid
and  strain within the solid.

The interaction forces and concentration forces satisfy the relations 
$\phi_s \mathbf  f_s  +\phi_f \mathbf f_s = 0$ and
$\phi_s \nabla \pi_s  + \phi_f \nabla \pi = 0$ \cite{lanir}.
The osmotic pressure is a function of the biomass fraction
$ \phi_f = \Pi (\phi_s)$ \cite{seminara}.
For isotropic solids with isotropic permeability the filtration force
\begin{eqnarray}
\mathbf f_f = - {1\over k_h} \mathbf q,
\label{filtrationforce}
\end{eqnarray}
where $k_h$ (hydraulic permeability) is a positive function
of $\phi_s$ \cite{lanir}.  Typically, $k_h(\phi_f)={\phi_f^2 \over \zeta}$, 
where  $\zeta$ is a friction parameter often set equal to
$\zeta= {\mu_f \over \xi(\phi_s)^2} >0$ and $\xi$ is the ``mesh size'' 
of the underlying biomass network \cite{seminara}.

Using the expressions for the stress tensors (\ref{stressfluid}) and 
(\ref{stresssolid}),  equations (\ref{motion}) become
\begin{eqnarray} 
{\rm div} \, \hat{\boldsymbol \sigma}_s  
+ \phi_s (-\nabla p + \nabla \pi_s ) + \phi_s \mathbf  f_s = 0,
\quad
\phi_f  (-\nabla p+\nabla \pi) + \phi_f  \mathbf  f_f = 0.
\label{motionsolidfluid}
\end{eqnarray}

Combining (\ref{motionsolidfluid}), (\ref{filtrationforce}), and (\ref{flux})
we obtain
\begin{eqnarray}
\mathbf q = - k_h \nabla (p - \pi) = \phi_f (\mathbf v_f -\mathbf v_s).
\label{darcy}
\end{eqnarray}
This is Darcy's law in the presence of concentration gradients.

Adding up equations (\ref{motionsolidfluid}), we find an equation relating 
solid displacements and pressure
\begin{eqnarray}
{\rm div} \, \hat{\boldsymbol \sigma}_s(\mathbf u_s)
- \nabla p = 0.
\label{motionsolid2}
\end{eqnarray}
At the biofilm boundary, the jumps in the total
stress vector  and  the chemical potential vanish:
\begin{eqnarray*}
(\hat{\boldsymbol \sigma}_s  - p \mathbf I) \mathbf n 
= \mathbf t_{ext},  \quad
p - \pi = p_{ext} - \pi_{f,ext}, 
\end{eqnarray*}
when applicable.

The solid velocity is then $\mathbf v_s ={\partial \mathbf u_s
\over \partial t}.$ These equations are complemented by
 (\ref{lcsolid}) and (\ref{balancemixture}), which now becomes
\begin{eqnarray}
{\rm div}(\mathbf v_s) = - {\rm div}(\mathbf q) =
{\rm div}(k_h \nabla (p - \pi)).
\label{balancemixture2}
\end{eqnarray}

\vskip 5mm

{\bf Acknowledgements.} 
This research has been partially supported by the FEDER /Ministerio
de Ciencia, Innovación y Universidades - Agencia Estatal de Investigación
grant PID2020-112796RB-C21.

\end{document}